\documentclass[10pt]{elsarticle}
 
\usepackage{hyperref}

\usepackage{calrsfs,amsmath,amssymb,graphicx,array,accents,figlatex,epsfig}
\usepackage{setspace,wrapfig,colortbl,marvosym, bbm, stmaryrd,centernot, color,upgreek}
\usepackage{amsfonts}
\usepackage{epstopdf}
\usepackage{algorithmic}
\usepackage{adjustbox}
\usepackage{multirow}

\ifpdf
  \DeclareGraphicsExtensions{.eps,.pdf,.png,.jpg}
\else
  \DeclareGraphicsExtensions{.eps}
\fi

\usepackage{enumitem}
\setlist[enumerate]{leftmargin=.5in}
\setlist[itemize]{leftmargin=.5in}

\newlength{\kaka}

\newcommand{\ahref}[2]{}

\newcommand{\beq}{\begin{equation}}
\newcommand{\eeq}{\end{equation}}
\newcommand{\lb}{\label}

\newcommand{\bea}{\begin{eqnarray}}
\newcommand{\eea}{\end{eqnarray}}
\newcommand{\bxr}{\begin{array}}
\newcommand{\exr}{\end{array}}

\newcommand\exs{\hspace*{0.4mm}}
\newcommand\xxs{\hspace*{0.2mm}}

\newcommand\nxs{\hspace*{-0.2mm}}

\newcommand{\norms}[1]{\parallel\! #1 \!\parallel}

\newcommand{\bn} {\boldsymbol{n}}

\newcommand{\bx} {\boldsymbol{x}}

\newcommand{\btau} {{\boldsymbol{\tau}}}

\newcommand{\bd} {\boldsymbol{d}}
\newcommand{\bI} {\boldsymbol{I}}

\newcommand{\sip} {\!\cdot\!}

\newcommand{\bzero}{\boldsymbol{0}}

\newcommand{\bu} {\boldsymbol{u}}

\newcommand{\bxi} {\boldsymbol{\xi}}

\begin{document}

\begin{frontmatter}

\title{Deep learning for full-field ultrasonic characterization}

\author{Yang Xu$^1$}
\author{Fatemeh Pourahmadian$^{1,2}$\corref{cor1}} 
\author{Jian Song$^1$}
\author{Conglin Wang$^3$}

\address{\vspace*{1.25mm}$^1$ Department of Civil, Environmental \& Architectural Engineering, University of Colorado Boulder, USA} 
\address{\vspace*{-1mm}$^2$ Department of Applied Mathematics, University of Colorado Boulder, USA}
\cortext[cor1]{Corresponding author: tel. 303-492-2027, email {\tt fatemeh.pourahmadian@colorado.edu}}
\address{\vspace*{-1mm}$^3$ Department of Physics, University of Colorado Boulder, USA}

\begin{abstract}
This study takes advantage of recent advances in machine learning to establish a physics-based data analytic platform for distributed reconstruction of mechanical properties in layered components from full waveform data. In this vein, two logics, namely the direct inversion and physics-informed neural networks (PINNs), are explored. The direct inversion entails three steps:~(i)~spectral denoising and differentiation of the full-field data,~(ii)~building appropriate neural maps to approximate the profile of unknown physical and regularization parameters on their respective domains, and~(iii) simultaneous training of the neural networks by minimizing the Tikhonov-regularized PDE loss using data from (i). PINNs furnish efficient surrogate models of complex systems with predictive capabilities via multitask learning where the field variables are modeled by neural maps endowed with (scaler or distributed) auxiliary parameters such as physical unknowns and loss function weights. PINNs are then trained by minimizing a measure of data misfit subject to the underlying physical laws as constraints. In this study, to facilitate learning from ultrasonic data, the PINNs loss adopts (a) wavenumber-dependent Sobolev norms to compute the data misfit, and (b) non-adaptive weights in a specific scaling framework to naturally balance the loss objectives by leveraging the  form of PDEs germane to elastic-wave propagation. Both paradigms are examined via synthetic and laboratory test data. In the latter case, the reconstructions are performed at multiple frequencies and the results are verified by a set of complementary experiments highlighting the importance of verification and validation in data-driven modeling.        
\end{abstract}

\end{frontmatter}

\section{Introduction}

Recent advances in laser-based ultrasonic testing has led to the emergence of dense spatiotemporal datasets which along with suitable data analytic solutions may lead to better understanding of the mechanics of complex materials and components. This includes learning of distributed mechanical properties from test data which is of interest in a wide spectrum of applications from medical diagnosis to additive manufacturing~\cite{liang2009acoustomotive, bal2014, garra2015elastography, bell2016, wei2021,chen2021learning,you2022}. This work makes use of recent progress in deep learning~\cite{bishop2006pattern, lecun2015deep} germane to direct and inverse problems in partial differential equations~\cite{cuomo2022scientific,wang2022(2),mccl2020,chen2018} to develop a systematic full-field inversion framework to recover the profile of pertinent physical quantities in layered components from laser ultrasonic measurements. The focus is on two paradigms, namely: the direct inversion and physics-informed neural networks (PINNs)~\cite{raissi2017physics, raissi2019physics, haghighat2021physics,henkes2022physics}. The direct inversion approach is in fact the authors' rendition of elastography method~\cite{muth1995,barb2004,baba2017} through the prism of deep learning. To this end, tools of signal processing are deployed to (a) denoise the experimental data, and (b) carefully compute the required field derivatives as per the governing equations. In parallel, the unknown distribution of PDE parameters in space-frequency are identified by neural networks which are then trained by minimizing the single-objective elastography loss. The learning process is stabilized via the Tikhonov regularization~\cite{tikh1995, kirsch2011} where the regularization parameter is defined in a distributed sense as a separate neural network which is simultaneously trained with the sought-for physical quantities. This unique exercise of learning the regularization field without a-priori estimates, thanks to neural networks, proved to be convenient, effective, and remarkably insightful in inversion of multi-fidelity experimental data.    
     
PINNs have recently come under the spotlight for offering efficient, yet predictive, models of complex PDE systems~\cite{cuomo2022scientific} that has so far been backed by rigorous theoretical justification within the context of linear elliptic and parabolic PDEs~\cite{shin2020}. Given the multitask nature of training for these networks and the existing challenges with modeling stiff and highly oscillatory PDEs~\cite{mccl2020,wang2022}, much of the most recent efforts has been focused on (a) adaptive gauging of the loss function~\cite{mccl2020,xian2022,bisc2021,son2022,zeng2022,yu2022,chen2018}, and (b) addressing the gradient pathologies~\cite{wang2022,chen2018} e.g., via learning rate annealing~\cite{wang2021} and customizing the network architecture~\cite{wang2022(2),jagtap2020,kim2022}. In this study, our initially austere implementations of PINNs using both synthetic and experimental waveforms led almost invariably to failure which further investigation attributed to the following impediments:~(a)~high-norm gradient fields due to large wavenumbers,~(b)~high-order governing PDEs in the case of laboratory experiments, and~(c)~imbalanced objectives in the loss function. These problems were further magnified by our attempts for distributed reconstruction of discontinuous PDE parameters -- in the case of laboratory experiments, from contaminated and non-smooth measurements. The following measures proved to be effective in addressing some of these challenges:~(i) training PINNs in a specific scaling framework where the dominant wavenumber is the reference length scale, (ii) using the wavenumber-dependent Sobolev norms in quantifying the data misfit, (iii) taking advantage of the inertia term in the governing PDEs to naturally balance the objectives in the loss function, and (iv) denoising of the experimental data prior to training.      

This paper is organized as follows. Section~\ref{Prelim} formulates the direct scattering problem related to the synthetic and laboratory experiments, and provides an overview of the data inversion logic. Section~\ref{Si} presents the computational implementation of direct inversion and PINNs to reconstruct the distribution of L\'{a}me parameters in homogeneous and heterogeneous models from in-plane displacement fields. Section~\ref{Le} provides a detailed account of laboratory experiments, scaling, signal processing, and inversion of antiplane particle velocity fields to recover the distribution of a physical parameter affiliated with flexural waves in thin plates. The reconstruction results are then verified by a set of complementary experiments. 

\vspace*{-2mm}
\section{Concept}\label{Prelim}

This section provides~(i)~a generic formalism for the direct scattering problem pertinent to the ensuing (synthetic and experimental) full-field characterizations, and~(ii)~data inversion logic. 

\vspace*{-1mm}
\subsection{Forward scattering problem}\label{FP}

Consider ultrasonic tests where the specimen $\Uppi \subset \mathbb{R}^d$, $d = 2,3$, is subject to (boundary or internal) excitation over the incident surface $S^\text{inc}\! \subset \overline{\Uppi}$ and the induced (particle displacement or velocity) field $\text{\bf{u}}\colon \overline{\Uppi} \times [0 \,\, T] \to \mathbb{R}^{N_\Lambda}\!$ ($N_\Lambda \leqslant d$) is captured over the observation surface $S^\text{obs}\! \subset \Uppi$ in a timeframe of length $T$. Here, $\Uppi$ is an open set whose closure is denoted by $\overline{\Uppi}$, and the sensing configuration is such that $\overline{S^\text{inc}\nxs} \cap \overline{S^\text{obs}\!} = \emptyset$. In this setting, the spectrum of observed waveforms $\hat{\text{\bf{u}}}\colon S^\text{obs}\nxs \times \Omega \to \mathbb{C}^{N_\Lambda}\!$ is governed by
\beq\lb{uf}
\Lambda[\hat{\text{\bf{u}}};{\boldsymbol{\vartheta}}](\bxi,\omega) ~=~ \bzero, \quad \hat{\text{\bf{u}}}~\colon\!\!\!=~\mathcal{F}[\text{\bf{u}}](\bxi,\omega), \quad \bxi \in S^\text{obs}\!, \, \omega \in \Omega, 
\eeq  
where $\Lambda$ of size $N_\Lambda \nxs\times 1$ designates a differential operator in frequency-space; $\mathcal{F}$ represents the temporal Fourier transform; ${\boldsymbol{\vartheta}}$ of dimension $N_\vartheta \times 1$ is the vector of relevant geometric and elastic parameters e.g.,~Lam{\'e} constants and mass density; $\bxi \in \mathbb{R}^d$ is the position vector; and $\omega>0$ is the frequency of wave motion within the specified bandwidth $\Omega$.     

\vspace*{-1mm}
\subsection{Dimensional platform}\label{DP}
All quantities in~\eqref{uf} are rendered dimensionless by identifying $\rho_\circ$, $\sigma_\circ$, and $\ell_\circ$ as the respective reference scales~\cite{Scaling2003} for mass density, elastic modulus, and length whose explicit values will be later specified. 

\vspace*{-1mm}
\subsection{Data inversion}\label{IP}

Given the full waveform data $\hat{\text{\bf{u}}}$ on $S^\text{obs}\nxs \times \Omega$, the goal is to identify the distribution of material properties over $S^\text{obs}$. For this purpose, two reconstruction paradigms based on neural networks are adopted in this study, namely:~(i) direct inversion, and~(ii) physics-based neural networks. Inspired by the elastography method~\cite{muth1995,barb2004}, quantities of interest in (i) are identified by neural maps over $S^\text{obs}\nxs \times \Omega$ that minimize a regularized measure of $\Lambda$ in~\eqref{uf}. The neural networks in (ii), however, are by design predictive maps of the waveform data (i.e.,~$\hat{\text{\bf{u}}}$) obtained by minimizing the data mismatch subject to~\eqref{uf} as a soft or hard constraint. In this setting, the unknown properties of $\Lambda$ may be recovered as distributed parameters of the (data) network during training via multitask optimization. In what follows, a detailed description of the deployed cost functions in (i) and (ii) is provided after a brief review of the affiliated networks.     

\subsubsection{Waveform and parameter networks}\label{ANN}

Laser-based ultrasonic experiments furnish a dense dataset on $S^\text{obs}\nxs \times \Omega$. Based on this, multilayer perceptrons (MLPs) owing to their dense range~\cite{horn1991} may be appropriate for approximating complex wavefields and distributed PDE parameters. Moreover, this architecture has proven successful in numerous applications within the PINN framework~\cite{raissi2019physics}. In this study, MLPs serve as both data and property maps where the input consists of discretized space and frequency coordinates $(\bxi_i,\omega_j)$, $i = 1,2, \ldots, N_\xi$, $j = 1,2, \ldots, N_\omega$, as well as distinct experimental parameters, e.g.,~the source location, distilled as one vector $\tau_k$ on domain $\mathcal{T}$ with $k = 1,2, \ldots, N_{\tau}$, while the output represents waveform data $\boldsymbol{\mathcal{D}}_{ijk} = [\mathfrak{R}{\hat{\text{\bf{u}}}}, \mathfrak{I}{\hat{\text{\bf{u}}}}](\bxi_i,\omega_j;\tau_k) \in \mathbb{R}^{N_\Lambda\!} \nxs\times \mathbb{R}^{N_\Lambda\!}$, and/or the sought-for mechanical properties $\boldsymbol{\mathcal{P}}_{ijn} = [\mathfrak{R}{\vartheta}_n, \mathfrak{I}{\vartheta}_n](\bxi_i,\omega_j) \in \mathbb{R} \nxs\times \mathbb{R}$, $n = 1,2, \ldots, N_\vartheta$. Note that following~\cite{chen2022}, the real $\mathfrak{R}$ and imaginary $\mathfrak{I}$ parts of~\eqref{uf} and every complex-valued variable are separated such that both direct and inverse problems are reformulated in terms of real-valued quantities. In this setting, each fully-connected MLP layer with $N_l$ neurons is associated with the forward map $\Upsilon_l\colon \mathbb{R}^{N_{l-1}} \to \mathbb{R}^{N_{l}}$,
\begin{equation}
\Upsilon_l(\boldsymbol{x}^{l-1}) ~=~ \text{tanh}(\boldsymbol{W}^l \boldsymbol{x}^{l-1} +\, \boldsymbol{b}^l), \quad \boldsymbol{x}^{l-1} \in\, \mathbb{R}^{N_{l-1}},
\end{equation}
where $\boldsymbol{W}^l \in \mathbb{R}^{N_l\times N_{l-1}}$ and $\boldsymbol{b}^l \in \mathbb{R}^{N_l}$ respectively denote the $l^\text{th}$ layer's weight and bias. Consecutive composition of $\Upsilon_l$ for $l = 1,2, \ldots, N_m$ builds the network map wherein $N_m$ designates the number of layers.    

\subsubsection{Direct inversion}\label{DIT}

\begin{figure}[!bp]
\vspace*{-3mm}
\center\includegraphics[width=0.64\linewidth]{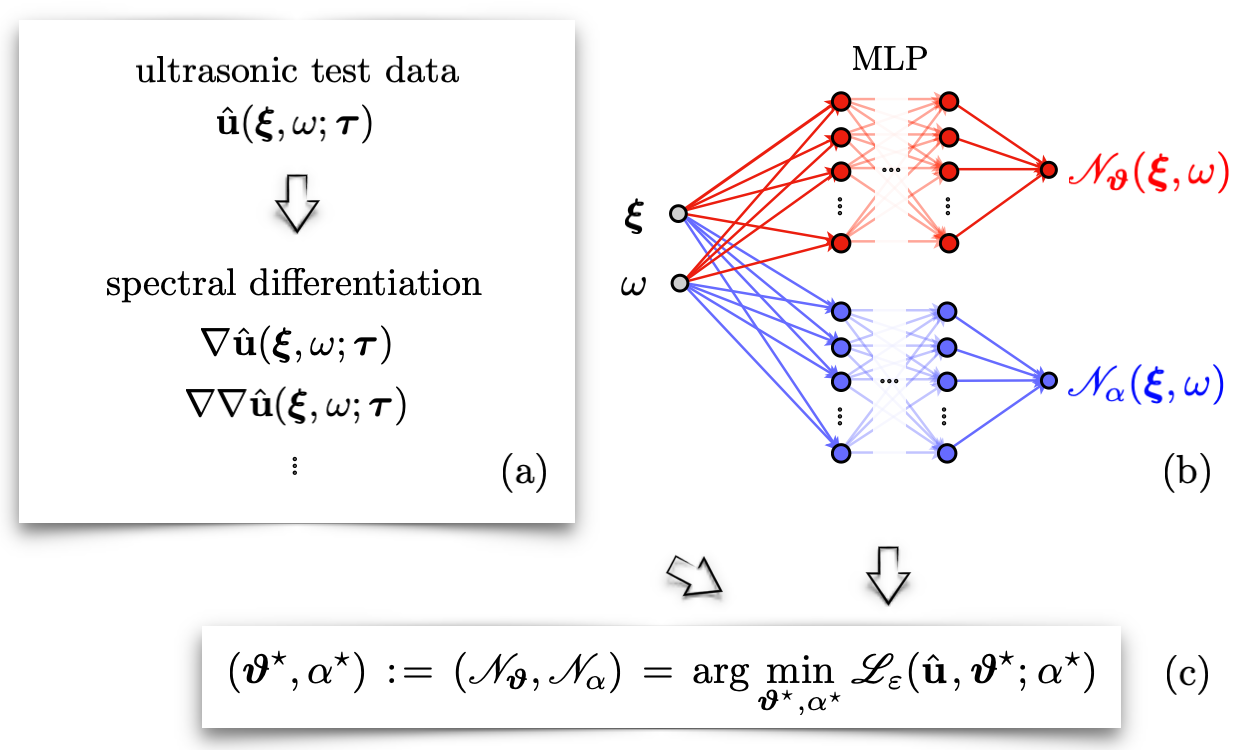} \vspace*{-2.5mm} 
\caption{Direct inversion:~(a) FFT-based spatial differentiation of the full-field data as per operator $\Lambda$,~(b)~MLP-based approximation of the unknown PDE and regularization parameters $({\boldsymbol{\vartheta}}, \alpha)$ on their respective domains, and (c) training the MLPs via minimizing the elastography loss $\mathcal{L}_\varepsilon$ according to~\eqref{Lels}.}
\label{DIDM7}
\vspace*{-5mm}
\end{figure} 

Logically driven by the elastography method, the direct inversion approach depicted in~Fig.~\ref{DIDM7} takes advantage of the leading-order physical principles underpinning the test data to recover the distribution of relevant physical quantities in space-frequency i.e.,~over the measurement domain. The ML-based direct inversion entails three steps:~(a) spectral denoising and differentiation of (n-differentiable) waveforms $\hat{\text{\bf{u}}}$ over $S^\text{obs}\nxs \times \Omega$ according to the (n-th order) governing PDEs in~\eqref{uf},~(b) building appropriate MLP maps to estimate the  profile of unknown physical parameters of the forward problem and regularization parameters of the inverse solution, and~(c) learning the MLPs through regularized fitting of data to the germane PDEs.       

Note that synthetic datasets -- generated via e.g.,~computer modeling or the method of manufactured solutions, may directly lend themselves to the fitting process in~(c) as they are typically smooth by virtue of numerical integration or analytical form of the postulated solution. Laboratory test data, however, are generally contaminated by noise and uncertainties, and thus, spectral differentiation is critical to achieve the smoothness requirements in~(c). The four-tier signal processing of experimental data follows closely that of~\cite[Section 3.1]{pour2018(2)} which for completeness is summarized here:~(1)~a band-pass filter consistent with the frequency spectrum of excitation is applied to the measured time signals at every receiver point,~(2) the obtained temporally smooth signals are then differentiated or integrated to obtain the pertinent field variables,~(3) spatial smoothing is implemented at every snapshot in time via application of median and moving average filters followed by computing the Fourier representation of the processed waveforms in space,~(4) the resulting smooth fields may be differentiated (analytically in the Fourier space) as many times as needed based on the underlying physical laws in preparation for the full-field reconstruction in step (c).          
It should be mentioned that the experimental data may feature intrinsic discontinuities e.g.,~due to material heterogeneities or contact interfaces. In this case, the spatial smoothing in (3) must be implemented in a piecewise manner after the geometric reconstruction of discontinuity surfaces in $S^\text{obs}$ which is quite straightforward thanks to the full-field measurements, see e.g.,~\cite[section 3.2]{pour2018(2)}. 

Next, the unknown PDE parameters ${\boldsymbol{\vartheta}}$ are approximated by a fully connected MLP network ${\boldsymbol{\vartheta}}^\star \colon \!\!\! = \mathcal{N}_{\boldsymbol{\vartheta}}(\bxi,\omega)$ as per Section~\ref{ANN}. The network is trained by minimizing the loss function    
\begin{equation}\lb{Lels}
\mathcal{L}_\varepsilon(\hat{\text{\bf{u}}}, {\boldsymbol{\vartheta}}^\star;\alpha) ~=~ \lVert \Lambda(\hat{\text{\bf{u}}};{\boldsymbol{\vartheta}}^\star) \rVert^2_{L^2(S^\text{obs}\nxs \times\xxs \Omega \xxs\times \mathcal{T})^{N_\Lambda\!}} \,+\,  \lVert \alpha \boldsymbol{\mathbbm{1}}_{\vartheta} \nxs\odot\nxs  {\boldsymbol{\vartheta}}^\star \rVert^2_{L^2(S^\text{obs}\nxs \times\xxs \Omega)^{N_\vartheta}},
\end{equation}
where $\boldsymbol{\mathbbm{1}}_{\vartheta}$ indicates an all-ones vector of dimension $N_{{\vartheta}} \times 1$, and $\odot$ designates the (element-wise) Hadamard product. Here, the PDE residual based on~\eqref{uf} is penalized by the norm of unknown parameters. Observe that the latter is a function of the weights and biases of the neural network which may help stabilize the MLP estimates during optimization. Such Tikhonov-type functionals are quite common in waveform tomography applications~\cite{liu2022, Fatemeh2017,cako2016} owing to their well-established regularizing properties~\cite{tikh1995,kirsch2011}. Within this framework, $\mathbb{R} \ni \alpha > 0$ is the regularization parameter which may be determined by three means, namely: (i) the Morozov discrepancy principle~\cite{moro2012, Kress1999}, (ii) its formulation as a (constant or distributed) parameter of the ${\boldsymbol{\vartheta}}^\star\!$ network which could then be learned during training, and (iii) its independent reconstruction as a separate MLP network ${{\alpha}}^\star \colon \!\!\! = \mathcal{N}_{\alpha}(\bxi,\omega)$ illustrated in~Fig.~\ref{DIDM7}~(b) that is simultaneously trained along with ${\boldsymbol{\vartheta}}^\star\nxs$ by minimizing~\eqref{Lels}. In this study, direct inversion is applied to synthetic and laboratory test data with both $\alpha = 0$ and $\alpha > 0$, based on (ii) and (iii). It was consistently observed that the regularization parameter $\alpha$ plays a key role in controlling the MLP estimates. This is particularly the case in situations where the field $\hat{\text{\bf{u}}}$ is strongly polarized or near-zero in certain neighborhoods which brings about instability i.e.,~very large estimates for ${\boldsymbol{\vartheta}}^\star\!$ in these areas. In light of this, all direct inversion results in this paper correspond to the case of $\alpha > 0$ identified by the MLP network ${{\alpha}}^\star$.

\subsubsection{Physics-informed neural networks}\label{PINN}

\begin{figure}[!tp]
\vspace*{-3mm}
\center\includegraphics[width=0.94\linewidth]{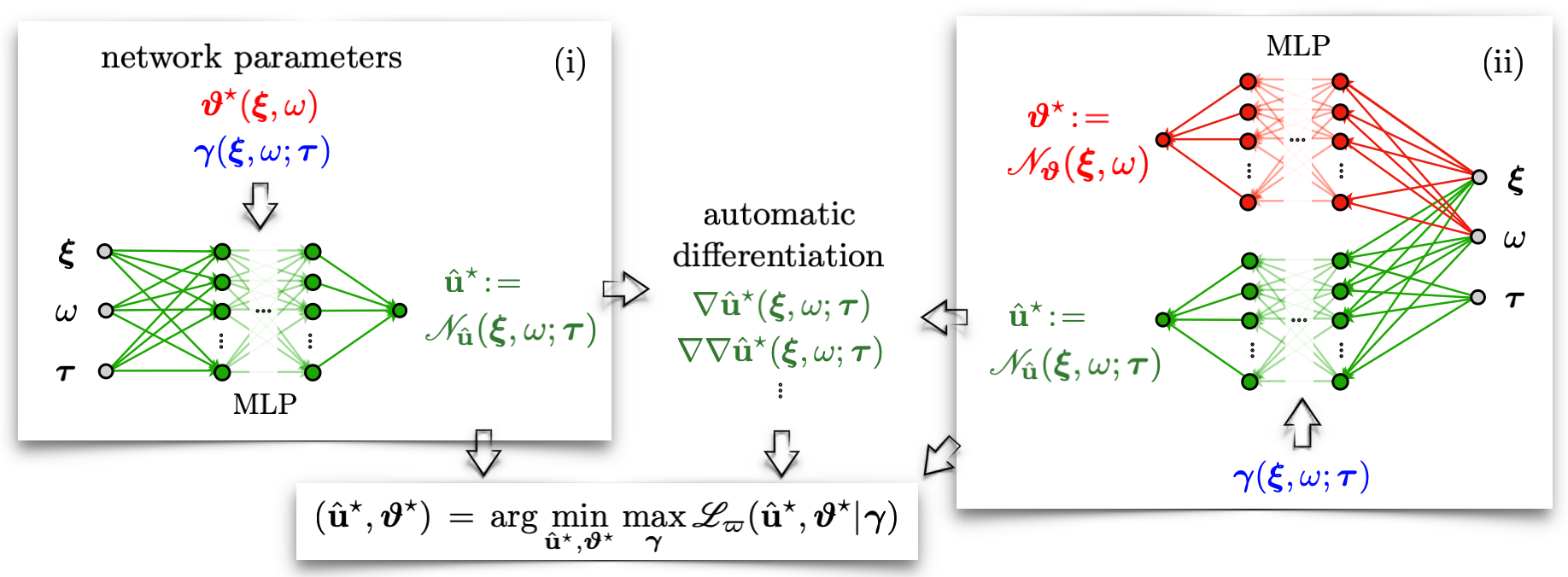} \vspace*{-3mm} 
\caption{Two logics for the physics-informed neural networks (PINNs) with distributed parameters:~(i)~the test data $\hat{\text{\bf{u}}}(\bxi,\omega;\boldsymbol{\tau})$ are modeled by a MLP map, while the unknown physical parameters ${\boldsymbol{\vartheta}}$ -- on $S^\text{obs}\nxs \times \Omega$, and the loss function weight ${{\gamma}}$ -- on $S^\text{obs}\nxs \times\xxs \Omega \xxs\times \mathcal{T}$, are defined as network parameters, and~(ii)~$\hat{\text{\bf{u}}}(\bxi,\omega;\boldsymbol{\tau})$ and ${\boldsymbol{\vartheta}}(\bxi,\omega)$ are identified by separate MLPs, while ${{\gamma}}$ is a parameter of $\mathcal{N}_{\hat{\text{\bf{u}}}}$. The MLP(s) in (i) and (ii) are then trained by minimizing $\mathcal{L}_{\varpi}$ of~\eqref{LPINNs} in the space of data and PDE parameters.}
\label{DIPINN4}
\vspace*{-7mm}
\end{figure}

By deploying the knowledge of underlying physics, PINNs~\cite{raissi2017physics,raissi2019physics} furnish efficient neural models of complex PDE systems with predictive capabilities. In this vein, a multitask learning process is devised according to Fig.~\ref{DIPINN4} where (a) the field variable $\hat{\text{\bf{u}}}$ -- i.e.,  measured data on $S^\text{obs}\nxs \times\xxs \Omega \xxs\times \mathcal{T}$, is modeled by the MLP map $\hat{\text{\bf{u}}}^\star \colon \!\!\! = \mathcal{N}_{\hat{\text{\bf{u}}}}(\bxi,\omega; \boldsymbol{\tau})$ endowed with the auxiliary parameter ${\gamma}(\bxi,\omega;\boldsymbol{\tau})$ related to the loss function~\eqref{LPINNs}, (b) the physical unknowns ${\boldsymbol{\vartheta}}$ could be defined either as parameters of $\hat{\text{\bf{u}}}^\star\!$ as in Fig.~\ref{DIPINN4}~(i), or as a separate MLP ${\boldsymbol{\vartheta}}^\star \colon \!\!\! = \mathcal{N}_{\boldsymbol{\vartheta}}(\bxi,\omega)$ as shown in Fig.~\ref{DIPINN4}~(ii), and~(c) learning the MLPs and affiliated parameters through minimizing a measure of data misfit subject to the governing PDEs as soft/hard constraints wherein the spatial derivatives of $\hat{\text{\bf{u}}}^\star\nxs$ are computed via automatic differentiation~\cite{pasz2017}. It should be mentioned that in this study all MLP networks are defined on (a subset of) $S^\text{obs}\nxs \times\xxs \Omega \xxs\times \mathcal{T}$ where $S^\text{obs} \nxs\cap \partial\Uppi= \emptyset$. Hence, the initial and boundary conditions -- which could be specified as additional constraints in the loss function~\cite{raissi2019physics}, are ignored. In this setting, the PINNs loss takes the form   
\begin{equation}\lb{LPINNs}
\mathcal{L}_\varpi(\hat{\text{\bf{u}}}^\star, {\boldsymbol{\vartheta}}^\star |\xxs {\gamma}) ~=~  \lVert \hat{\text{\bf{u}}} \exs - \exs \hat{\text{\bf{u}}}^\star \nxs \rVert^2_{\mathfrak{N}(S^\text{obs}\nxs \times\xxs \Omega \xxs\times \mathcal{T})^{N_\Lambda\!}} \,\exs+\,  \lVert \gamma  \boldsymbol{\mathbbm{1}}_{\Lambda} \nxs\odot\nxs \Lambda(\hat{\text{\bf{u}}}^\star;{\boldsymbol{\vartheta}}^\star) \rVert^2_{L^2(S^\text{obs}\nxs \times\xxs \Omega \xxs\times \mathcal{T})^{N_\Lambda\!}}, \,\,\, \mathfrak{N} \,=\, L^2, \widehat{H}^\iota, \,\,\, \iota \leqslant \text{n},
\end{equation}
where $\boldsymbol{\mathbbm{1}}_{\Lambda}$ is a $N_{{\Lambda}}\! \times 1$ vector of ones; $\text{n}$ is the order of $\Lambda$, and $\widehat{H}^\iota$ denotes the adaptive $H^\iota$ norm defined by

\beq\lb{Hsk}
\| \cdot \|_{\widehat{H}^\iota}:=\sqrt{\sum_{1 \xxs\leqslant\xxs |\boldsymbol{e}| \xxs\leqslant\exs \iota} \!\!\!\! \gamma^{\boldsymbol{e}} \exs \lVert \nabla^{\boldsymbol{e}}(\xxs\cdot\xxs) \rVert_{L^2}^2+\lVert \exs\cdot\exs \rVert^2_{L^2}}, \quad \nabla^{\boldsymbol{e}} = \frac{\partial^{|\boldsymbol{e}|}}{\partial \xi_1^{e_1} \partial \xi_2^{e_2} \sip\exs\sip\exs\sip \xxs\exs \partial \xi_d^{e_d}}, \quad  |\boldsymbol{e}| \,\colon\!\!\!= \sum_{\text{i} = 1}^d e_\text{i}.
\eeq

Here, $\boldsymbol{e} \colon\!\!\!= \lbrace e_1, e_2, \ldots e_d \rbrace$ is a vector of integers $e_\text{i} \geqslant 0$. Provided that $\forall\boldsymbol{e}, \,\gamma^{\boldsymbol{e}} = 1$, then ${\widehat{H}^\iota}$ is by definition equal to ${{H}^\iota}$~\cite{brez2010}. Note however that at high wavenumbers, ${{H}^\iota}$ is dominated by the highest derivatives $\nabla^{\boldsymbol{e}}\hat{\text{\bf{u}}}^\star$, $|\boldsymbol{e}| = \iota$, which may complicate (or even lead to the failure of) the training process due to uncontrolled error amplification by automatic differentiation particularly in earlier epochs. This issue may be addressed through proper weighting of derivatives in~\eqref{Hsk}. In light of the frequency-dependent Sobolev norms in~\cite{ha2003,liu2022}, one potential strategy is to adopt the wavenumber-dependent weights as the following
\beq\nonumber
\gamma^{\boldsymbol{e}} = \left(\frac{1}{\kappa_1^{e_1} \kappa_2^{e_2} \sip\exs\sip\exs\sip \xxs\exs \kappa_d^{e_d}}\right)^2,  \quad 1 \xxs\leqslant\xxs |\boldsymbol{e}| \xxs\leqslant\exs \iota, 
\eeq
wherein $\kappa_{\text{i}}$ is a measure of wavenumber along $\xi_{\text{i}}$ for ${\text{i}} =  1, \ldots, d$. In this setting, the weighted norms of derivatives in~\eqref{Hsk} remain approximately within the same order as the $L^2$ norm of data misfit. Another way to automatically achieve the latter is to set the reference scale $\ell_\circ$ such that $\kappa_{\text{i}} \!\sim\! 1$. Note that the ${\widehat{H}^\iota}$ norms directly inform the PINNs about the ``expected" field derivatives -- while preventing their uncontrolled magnification. This may help stabilize the learning process as such derivatives are intrinsically involved in the PINNs loss via $\Lambda(\hat{\text{\bf{u}}}^\star;{\boldsymbol{\vartheta}}^\star)$. It should be mentioned that when $\mathfrak{N} = \widehat{H}^\iota$ in~\eqref{LPINNs}, the ``true" estimates for derivatives $\nabla^{\boldsymbol{e}}\hat{\text{\bf{u}}}$ may be obtained via spectral differentiation as per Section~\ref{DIT}.  

The Lagrange multiplier~\cite{rockafellar1993lagrange,everett1963generalized} $\gamma(\bxi,\omega;\boldsymbol{\tau})$ in~\eqref{LPINNs} is critical for balancing the loss components during training. Its optimal value, however, highly depends on (a) the nature of $\Lambda$~\cite{mccl2020}, and (b) the distribution of unknown parameters $\boldsymbol{\vartheta}$. It should be mentioned that setting $\gamma = 1$ led to failure in almost all of the synthetic and experimental implementations of PINNs in this study. Gauging of loss function weights has been the subject of extensive recent studies~\cite{mccl2020,xian2022,liu2021,bisc2021,son2022,zeng2022}. One systematic approach is the adaptive SA-PINNs~\cite{mccl2020} where the multiplier $\gamma(\bxi,\omega;\boldsymbol{\tau})$ is a distributed parameter of $\hat{\text{\bf{u}}}^\star$ whose value is updated in each epoch according to a minimax weighting paradigm. Within this framework, the data (and parameter) networks are trained by minimizing $\mathcal{L}_\varpi$ with respect to $\hat{\text{\bf{u}}}^\star$ and ${\boldsymbol{\vartheta}}^\star$, while maximizing the loss with respect to $\gamma$ as shown in Fig.~\ref{DIPINN4}. 

Depending on the primary objective for PINNs, one may choose nonadaptive or adaptive weighting. More specifically, if the purpose is high-fidelity \emph{forward} modeling via neural networks where $\boldsymbol{\vartheta}$ is known a-priori and PINNs are intended to serve as predictive surrogate models of $\Lambda$, then ideas rooted in constrained optimization e.g.,~minimax weighting is theoretically sound. However, if the \emph{inverse} solution i.e., identification of $\boldsymbol{\vartheta}(\bxi,\omega)$ from ``real-world" or laboratory test data is the main goal particularly in a situation where any assumption on the smoothness of $\boldsymbol{\vartheta}$ and/or applicability of $\Lambda$ may be (at least locally) violated e.g., due to unknown material heterogeneities or interfacial discontinuities, then trying to enforce $\Lambda$ everywhere on $S^\text{obs}\nxs \times\xxs \Omega \xxs\times \mathcal{T}$ (via point-wise adaptive weighting) may lead to instability and failure of data inversion. In such cases, nonadaptive weighting may be more appropriate. In light of this, in what follows, $\gamma$ is a non-adaptive weight specified by taking advantage of the PDE structure to naturally balance the loss objectives.      

  \vspace*{-2mm} 
\section{Synthetic implementation}\label{Si}

\begin{figure}[!bp]
\vspace*{-4mm} 
\center\includegraphics[width=0.64\linewidth]{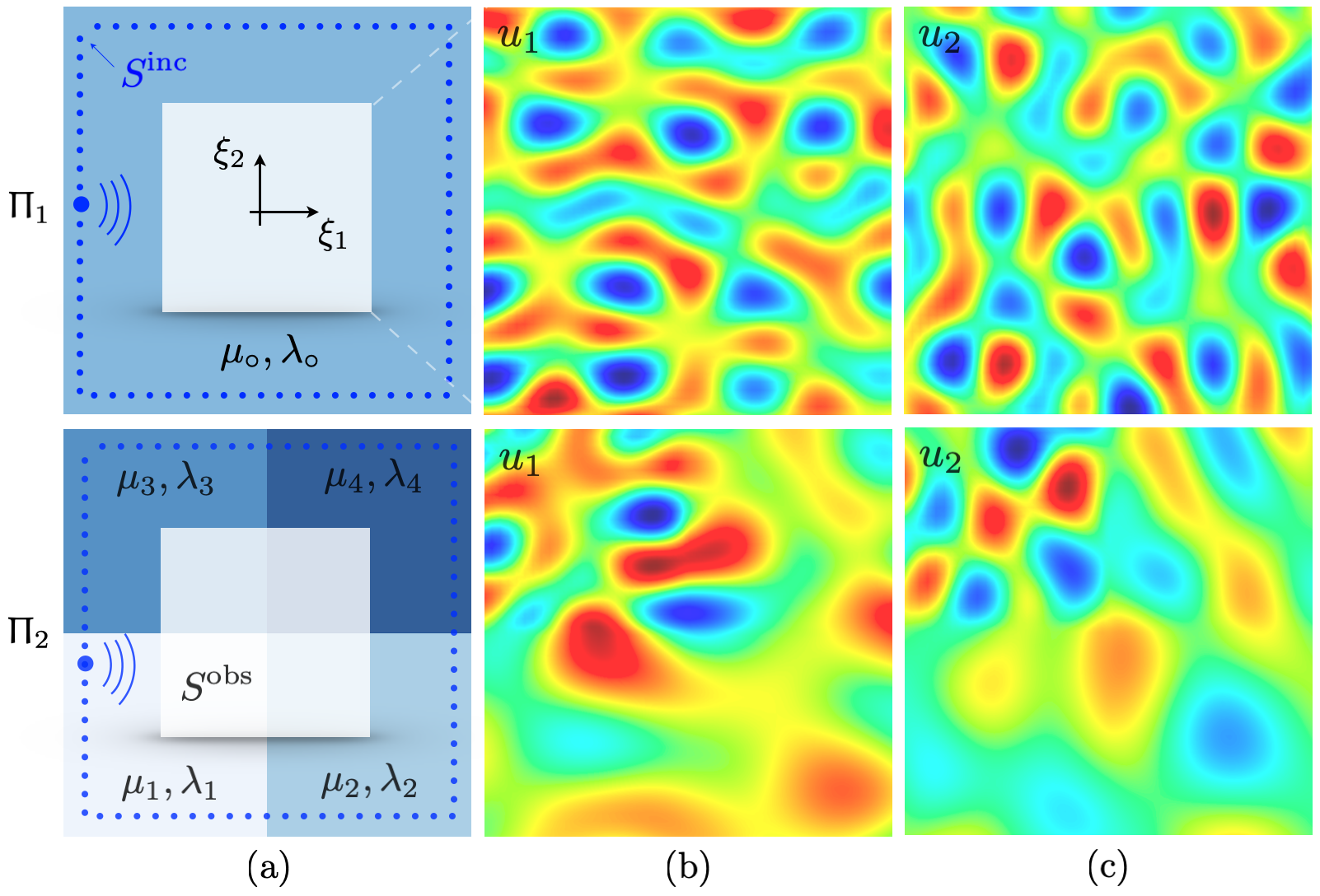} \vspace*{-3mm} 
\caption{synthetic experiments simulating plane-strain wave motion in homogeneous (top-left) and heterogeneous (bottom-left) specimens:~(a) testing configuration where the model is harmonically excited at frequency $\omega$ by a point source on $S^{\text{inc}}$, and the induced displacement field $\bu$ is computed over $S^{\text{obs}}$ along $\xi_1$ and~$\xi_2$ as shown in~(b) and (c), respectively.}\label{SE1}
\end{figure} 

Full-field characterization via the direct inversion and physics-informed neural networks are examined through a set of numerical experiments. The waveform data in this section are generated via a FreeFem++~\cite{Hech2012} code developed as part of~\cite{pour2022(1)}.

\vspace*{-2mm}
\subsection{Problem statement}\label{PS-synt}

Plane-strain wave motion in two linear, elastic, piecewise homogeneous, and isotropic samples is modeled according to~Fig.~\ref{SE1}~(a). On denoting the frequency of excitation by $\omega$, let $\ell_r = \frac{2\pi}{\omega} \sqrt{\mu_r/\rho_r}$, $\rho_r = 1$, and $\mu_r = 1$ be the reference scales for length, mass density, and stress, respectively. In this framework, both specimens are of size $16$ $\!\times\!\nxs$ $16$ and uniform density $\rho = 1$. The first sample $\Uppi_1 \subset \mathbb{R}^2$ is characterized by the constant Lam\'{e} parameters $\mu_\circ = 1$ and $\lambda_\circ = 0.47$, while the second sample $\Uppi_2 \subset \mathbb{R}^2$ is comprised of four perfectly bonded homogenous components $\Uppi_{2_j}\!$ of $\mu_j = j$ and $\lambda_j = 2j/3$, $j =  \lbrace 1, 2, 3, 4 \rbrace$ such that $\overline{\Uppi_2} = \bigcup_{j=1}^4 \overline{\Uppi_{2_j}}$. Accordingly, the shear and compressional wave speeds read $\textrm{c}_{\textrm{s}}^\circ = 1$, $\textrm{c}_{\textrm{p}}^\circ = 1.57$ in $\Uppi_1$, and $\textrm{c}_{\textrm{s}}^j = \sqrt{j}$, $\textrm{c}_{\textrm{p}}^j = 1.63\sqrt{j}$ in $\Uppi_{2_j}$. Every numerical experiment entails an in-plane harmonic excitation at $\omega = 3.91$ via a point source on $S^\text{inc}$ (the perimeter of a $14$ $\!\times\!\nxs$ $14$ square centered at the origin). The resulting displacement field $\bu^\upalpha = (u^\upalpha_1, u^\upalpha_2)$, $\upalpha = 1,2$, is then computed in $\Uppi_\upalpha$ over $S^\text{obs}$ (a concentric square of dimension $8$ $\!\times\!\nxs$ $8$) such that
\beq\lb{I-P} 
\begin{aligned}
&\mu_\upalpha \Delta \boldsymbol{u}^\upalpha(\bxi) \,+\, (\lambda_\upalpha\nxs+\mu_\upalpha)\nabla\nabla\nxs\cdot \boldsymbol{u}^\upalpha(\bxi) \,+\, \rho \xxs \omega^2 \boldsymbol{u}^\upalpha(\bxi)  ~=~ \delta(\bxi-\bx) \bd,  \quad & \bxi \in {\Uppi}_\upalpha, \bx \in S^\text{inc}, \\*[0.1mm]
&\big[ \lambda_\upalpha \nabla\nxs\cdot \boldsymbol{u}^\upalpha(\bxi) \bI_2  \,+\, 2\mu_\upalpha \nabla_\text{\tiny{sym}}\boldsymbol{u}^\upalpha(\bxi) \big]  \nxs\cdot \bn(\bxi) ~=~\bzero,  \quad & \bxi \in \partial {\Uppi}_\upalpha, \\*[0.0mm]
\end{aligned}   
\eeq     
where $\bx$ and $\bd$ respectively indicate the source location and polarization vector; $\bn$ is the unit outward normal to the specimen's exterior, and
\beq\label{par}\nonumber
\!\left\{\begin{array}{l}
\begin{aligned}
&\!\!  \mu_\alpha\,=\,\mu_\circ, \,\,\lambda_\alpha\,=\,\lambda_\circ, \quad   & \alpha = 1   \!\!\! \\*[0.1mm] 
& \!\! \mu_\alpha\,=\,\mu_j, \,\, \lambda_\alpha\,=\,\lambda_j, \quad   & \alpha = 2 \,\wedge\, \bxi \in \Uppi_{2_{j \in  \lbrace 1, 2, 3, 4 \rbrace}}   \!\!\! 
\end{aligned}
\end{array}\right..
\eeq

When $\upalpha = 2$, the first of~\eqref{I-P} should be understood as a shorthand for the set of four governing equations over $\Uppi_{2_j}$, $j =  \lbrace 1, 2, 3, 4 \rbrace$, supplemented by the continuity conditions for displacement and traction across $\partial \Uppi_{2_j} \!\!\setminus\! \partial \Uppi_{2}$ as applicable. 

In this setting, the generic form~\eqref{uf} may be identified as the following
\beq\label{map-synt}
\begin{aligned}
&\!\!  \Lambda ~=~\Lambda_\upalpha ~\colon\!\!\!=~ \mu_\upalpha \Delta  \,+\, (\lambda_\upalpha\nxs+\mu_\upalpha)\nabla\nabla\nxs\cdot  \,+\,\, \rho \xxs \omega^2 \bI_{2},  & \quad \upalpha ~=~ 1,2,  \\*[0.45mm] 
& \!\! \hat{\text{\bf{u}}} ~=~ \boldsymbol{u}^\upalpha(\bxi,\omega;\btau), \quad \boldsymbol{\vartheta} ~=~ [\xxs\mu_\upalpha,\lambda_\upalpha](\bxi,\omega),  & \quad \bxi \in S^\text{obs}\!, \, \omega \in \Omega, \btau \in \mathcal{T},
\end{aligned}
\eeq
wherein $\bI_{2}\!$ is the second-order identity tensor; $\btau = (\bx,\bd) \in S^{\text{inc}} \nxs\times\nxs \mathcal{B}_1 = \mathcal{T}$ with $\mathcal{B}_1$ denoting the unit circle of polarization directions. Note that $\rho$ is treated here as a \emph{known} parameter.

In the numerical experiments, $S^\text{inc}$ (\emph{resp}.~$S^\text{obs}$) is discretized by a uniform grid of 32 (\emph{resp}.~50 $\!\times\!$ 50) points, while $\Omega$ and $\mathcal{B}_1$ are respectively sampled at $\omega = 3.91$ and $\bd = (1,0)$.  

All inversions in this study are implemented within the PyTorch framework~\cite{PyT2019}. 

\subsection{Direct inversion}\label{DI-synt}

The three-tier logic of Section~\ref{DIT} is employed to reconstruct the distribution of $\mu_\upalpha$ and $\lambda_\upalpha$, $\upalpha = 1,2$, over $S^\text{obs}$, entailing:~(a)~spectral differentiation of the displacement field $\boldsymbol{u}^\upalpha$ in order to compute $\Delta \boldsymbol{u}^\upalpha$ and $\nabla\nabla\nxs\cdot \boldsymbol{u}^\upalpha$ as per~\eqref{I-P},~(b)~construction of three positive-definite MLP networks $\mu^\star$, $\lambda^\star$, and $\alpha^\star$; each of which is comprised of one hidden layer of 64 neurons, and~(c) training the MLPs by minimizing $\mathcal{L}_\varepsilon$ as in~\eqref{Lels} and~\eqref{map-synt} by way of the ADAM algorithm~\cite{kingma2014}. To avoid near-boundary errors affiliated with the one-sided FFT differentiation in $\Delta \boldsymbol{u}^\upalpha$ and $\nabla\nabla\nxs\cdot \boldsymbol{u}^\upalpha$, a concentric $40 \times 40$ subset of collocation points sampling $S^\text{obs}$ is deployed for training purposes. It should also be mentioned that in the heterogeneous case, i.e., $\upalpha = 2$, the discontinuity of derivatives across $\partial \Uppi_{2_{j \in  \lbrace 1, 2, 3, 4 \rbrace}}\!$ calls for piecewise spectral differentiation. According to Section~\ref{ANN}, the input to $\mathcal{P}^\star = \mathcal{N}_{\mathcal{P}}(\bxi,\omega)$, $\mathcal{P} = \mu, \lambda$, and $\alpha^\star = \mathcal{N}_{\alpha}(\bxi,\omega)$ is of size $N_\xi N_\tau \times N_\omega = 1600 N_s \times 1$ where $N_s \leqslant 32$ is the number of simulations i.e., source locations used to generate distinct waveforms for training. In this setting, since the physical quantities of interest are independent of $\btau$, the real-valued output of MLPs is of dimension $1600 \times 1$ furnishing a local estimate of the L\'{a}me and regularization parameters at the specified sampling points on $S^\text{obs}$. Each epoch makes use of the full dataset and the learning rate is $0.005$. 

In this work, the reconstruction error is measured in terms of the normal misfit               
 \begin{equation}\label{EM}
\Xi(\text{q}^\star) ~=~ \frac{\norms{\text{q}^\star\exs-\,\text{q}\xxs}_{L^2}}{\norms{\text{q}\xxs}_{L^\infty}},
\end{equation}
where $\text{q}^\star$ is an MLP estimate for a quantity with the ``true" value $\text{q}$.

\begin{figure}[!h]
\center\includegraphics[width=0.9\linewidth]{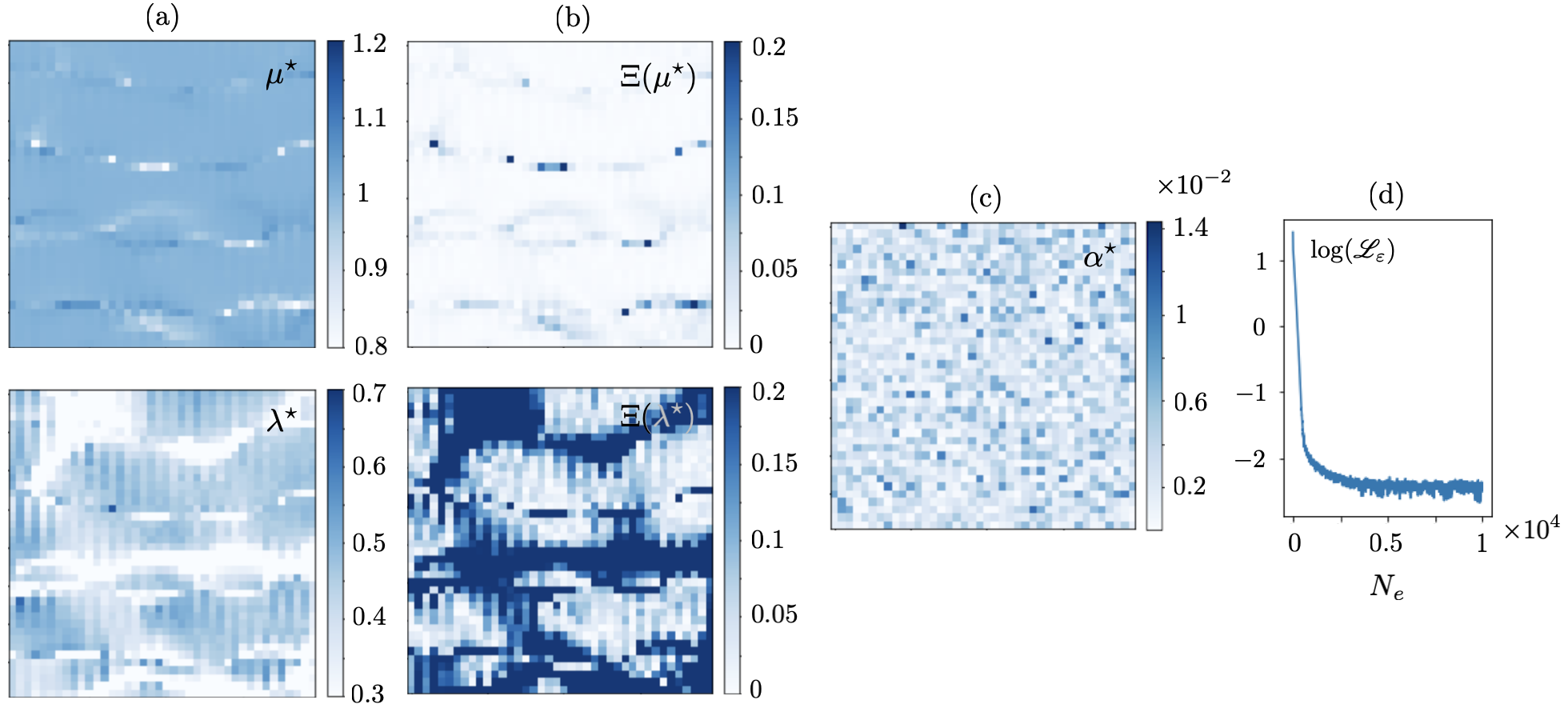} \vspace*{-2.5mm} 
\caption{Direct inversion of the L\'{a}me parameters in $\Uppi_1$ using noiseless data from a single source:~(a)~MLP-predicted distributions $\mu^\star\nxs$ and $\lambda^\star\nxs$,~(b)~reconstruction error~\eqref{EM} with respect to the true values $\mu_\circ = 1$ and $\lambda_\circ = 0.47$,~(c) MLP-recovered distribution of the regularization parameter $\alpha^\star\nxs$, and~(d) loss function $\mathcal{L}_\varepsilon$ \emph{vs}.~the number of epochs $N_e$ in the $\log = \log_{10}$ scale.}
\label{DR15}
\end{figure}

\begin{figure}[!h]
\center\includegraphics[width=0.9\linewidth]{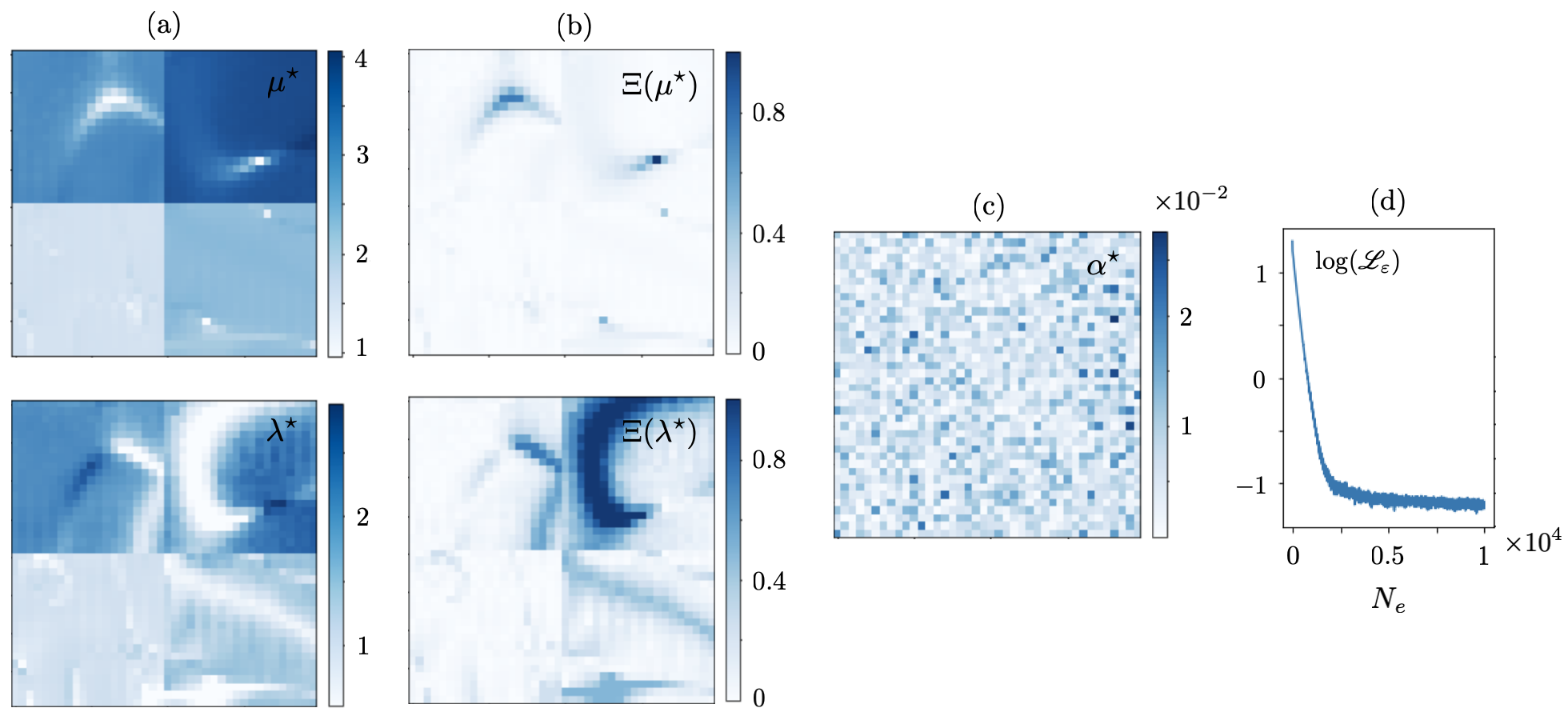} \vspace*{-2.5mm} 
\caption{Direct inversion of the L\'{a}me parameters in $\Uppi_2$ using noiseless data from a single source:~(a)~MLP-predicted distributions $\mu^\star\nxs$ and $\lambda^\star\nxs$,~(b)~reconstruction error~\eqref{EM} with respect to the true values $\mu_j = j$ and $\lambda_j = 2j/3$, $j =  \lbrace 1, 2, 3, 4 \rbrace$,~(c) MLP-recovered regularization parameter $\alpha^\star\nxs$, and~(d) loss function $\mathcal{L}_\varepsilon$ \emph{vs}.~the number of epochs $N_e$.}
\label{DR18}
\end{figure}

Let $S^\text{inc}$ be sampled at one point i.e.,~$N_s =1$ so that a single forward simulation in $\Uppi_\upalpha$, $\upalpha = 1,2$, generates the training dataset. The resulting reconstructions are shown in Figs.~\ref{DR15} and~\ref{DR18}. It is evident from both figures that the single-source reconstruction fails at the loci of near-zero displacement which may explain the relatively high values of the recovered regularization parameter $\alpha^\star$. Table~\ref{Num1} details the true values as well as mean and standard deviation of the reconstructed L\'{a}me distributions $\boldsymbol{\vartheta}^\star\nxs = (\mu^\star, \lambda^\star)$ in $\Uppi_1$ (\emph{resp}.~$\Uppi_{2_j}$ for $j = 1, 2, 3, 4$) according to Fig.~\ref{DR15} (\emph{resp}.~Fig.~\ref{DR18}).  

\begin{figure}[!h]
\center\includegraphics[width=0.9\linewidth]{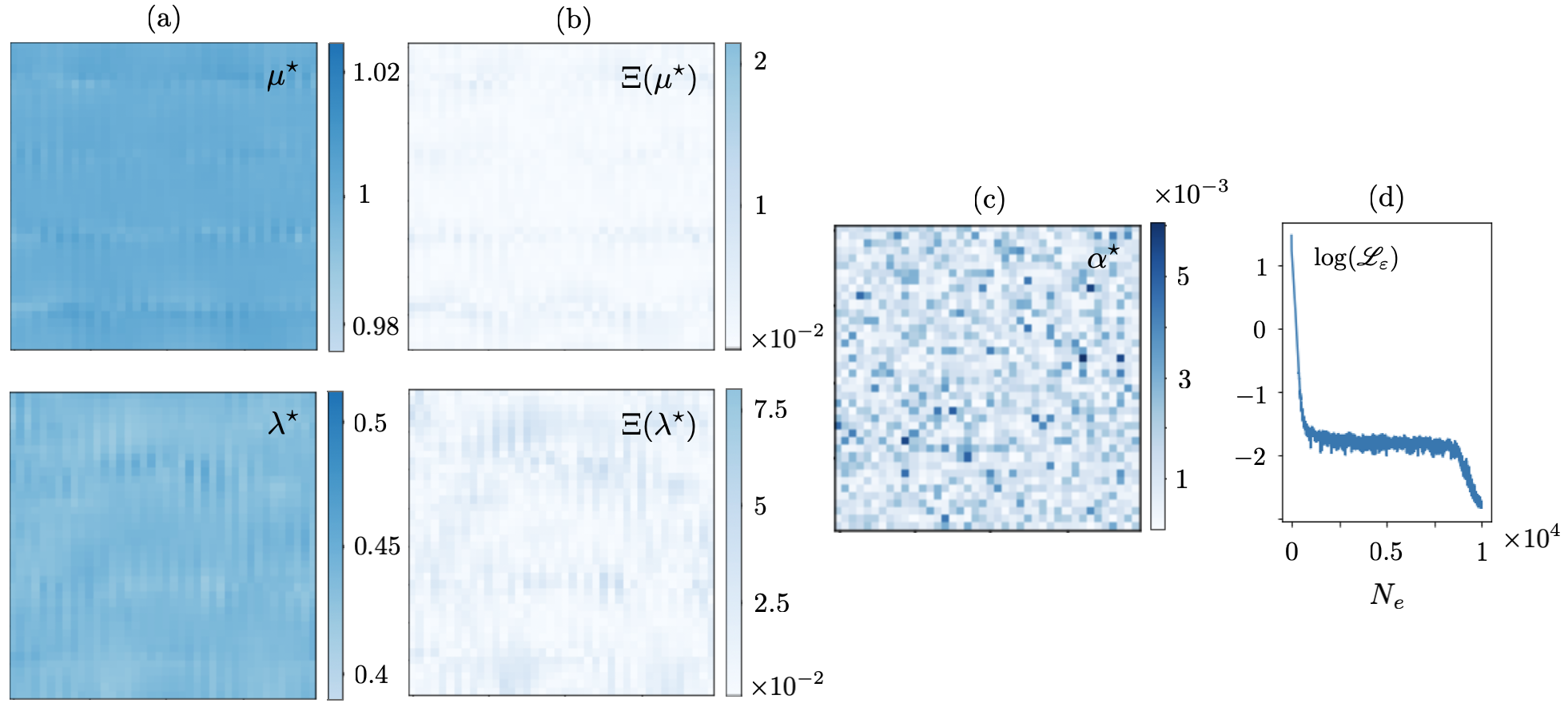} \vspace*{-2mm} 
\caption{Direct inversion of the L\'{a}me parameters in $\Uppi_1$ using noiseless data from five distinct simulations:~(a)~MLP-predicted distributions $\mu^\star\nxs$ and $\lambda^\star\nxs$,~(b)~reconstruction error~\eqref{EM} with respect to the true values $\mu_\circ = 1$ and $\lambda_\circ = 0.47$,~(c) MLP-recovered regularization parameter $\alpha^\star\nxs$, and~(d) loss function $\mathcal{L}_\varepsilon$ \emph{vs}.~the number of epochs $N_e$.}
\label{DR4}
\end{figure}

\begin{figure}[!h]
\center\includegraphics[width=0.9\linewidth]{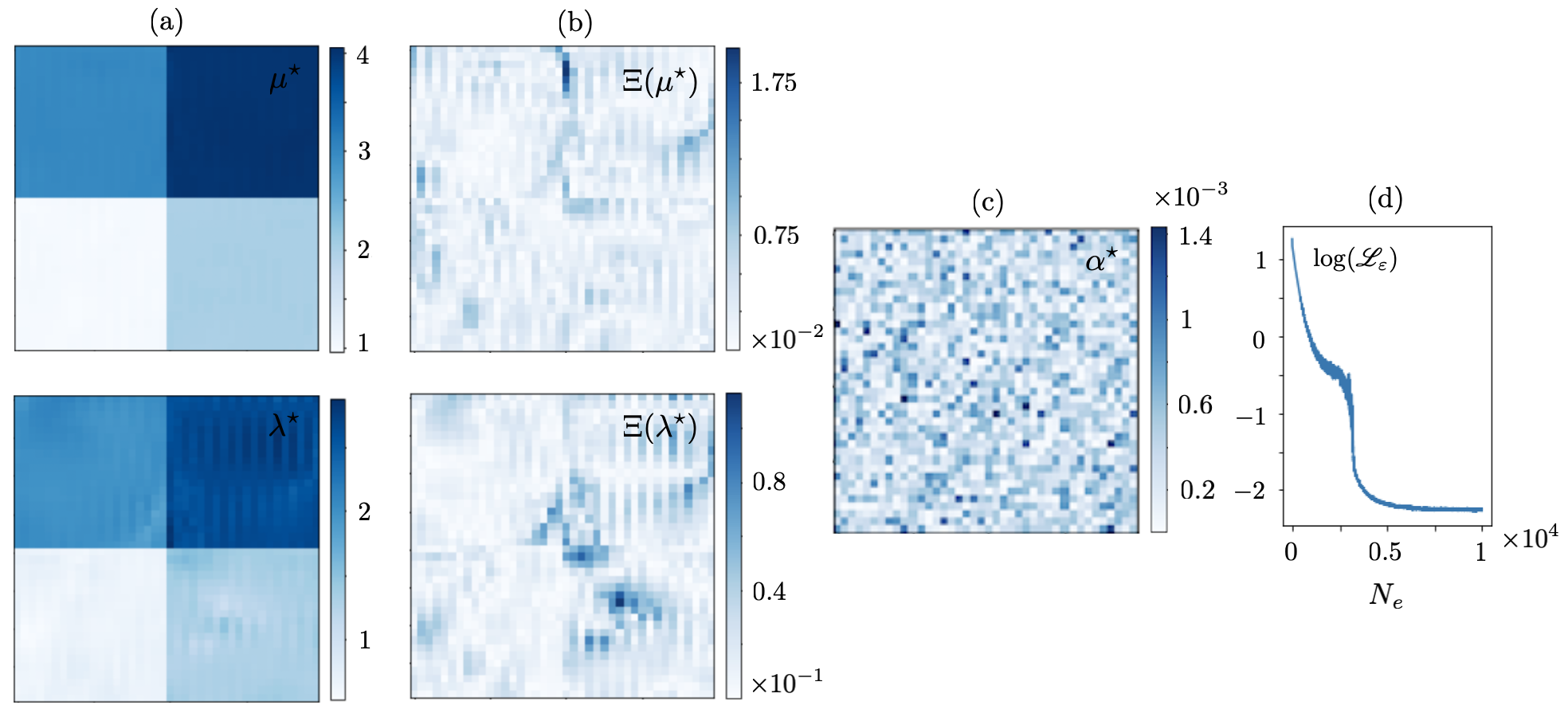} \vspace*{-2mm} 
\caption{Direct inversion of the L\'{a}me parameters in $\Uppi_2$ using five noiseless datasets:~(a)~MLP-predicted distributions $\mu^\star\nxs$ and $\lambda^\star\nxs$,~(b)~reconstruction error~\eqref{EM} with respect to the true values $\mu_j = j$ and $\lambda_j = 2j/3$, $j =  \lbrace 1, 2, 3, 4 \rbrace$,~(c) MLP-recovered regularization parameter $\alpha^\star\nxs$, and~(d) loss function $\mathcal{L}_\varepsilon$ \emph{vs}.~the number of epochs $N_e$.}
\label{DR6}
\end{figure}

This problem may be addressed by enriching the training dataset e.g., via increasing $N_s$. Figs.~\ref{DR4} and~\ref{DR6} illustrate the reconstruction results when $S^\text{inc}$ is sampled at $N_s =5$ source points. The mean and standard deviation of the reconstructed distributions are provided in Table~\ref{Num2}. It is worth noting that in this case the identified regularization parameter $\alpha^\star\nxs$ assumes much smaller values -- compared to that of Figs.~\ref{DR15} and~\ref{DR18}. This is closer to the scale of computational errors in the forward simulations.  

\begin{figure}[!h]
\center\includegraphics[width=0.9\linewidth]{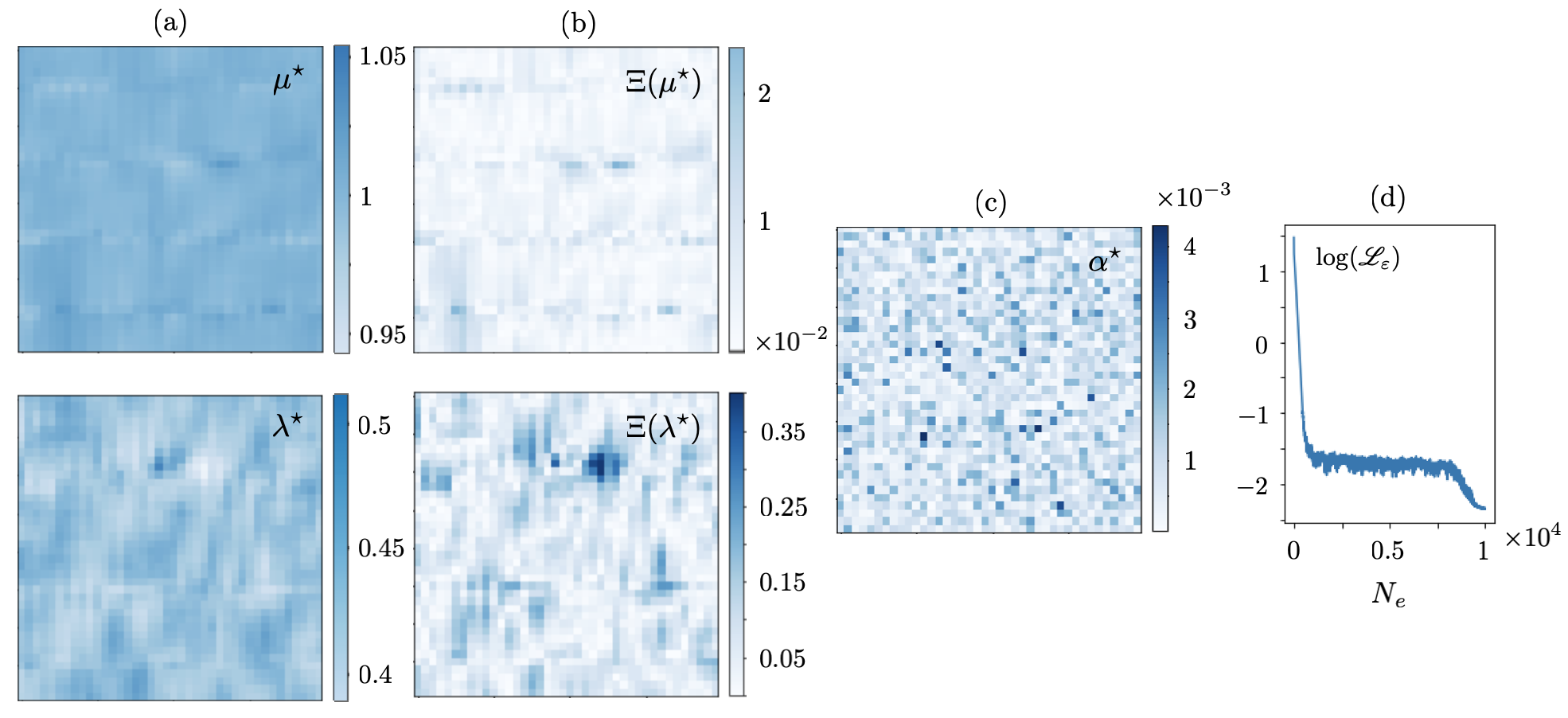} \vspace*{-2mm} 
\caption{Direct inversion of the L\'{a}me parameters in $\Uppi_1$ using five datasets perturbed with 5\% white noise:~(a)~MLP-predicted distributions $\mu^\star\nxs$ and $\lambda^\star\nxs$,~(b)~reconstruction error~\eqref{EM} with respect to the true values $\mu_\circ = 1$ and $\lambda_\circ = 0.47$,~(c) MLP-recovered regularization parameter $\alpha^\star\nxs$, and~(d) loss function $\mathcal{L}_\varepsilon$ \emph{vs}.~the number of epochs $N_e$.}
\label{DR5}
\end{figure}

\begin{figure}[!h]
\center\includegraphics[width=0.9\linewidth]{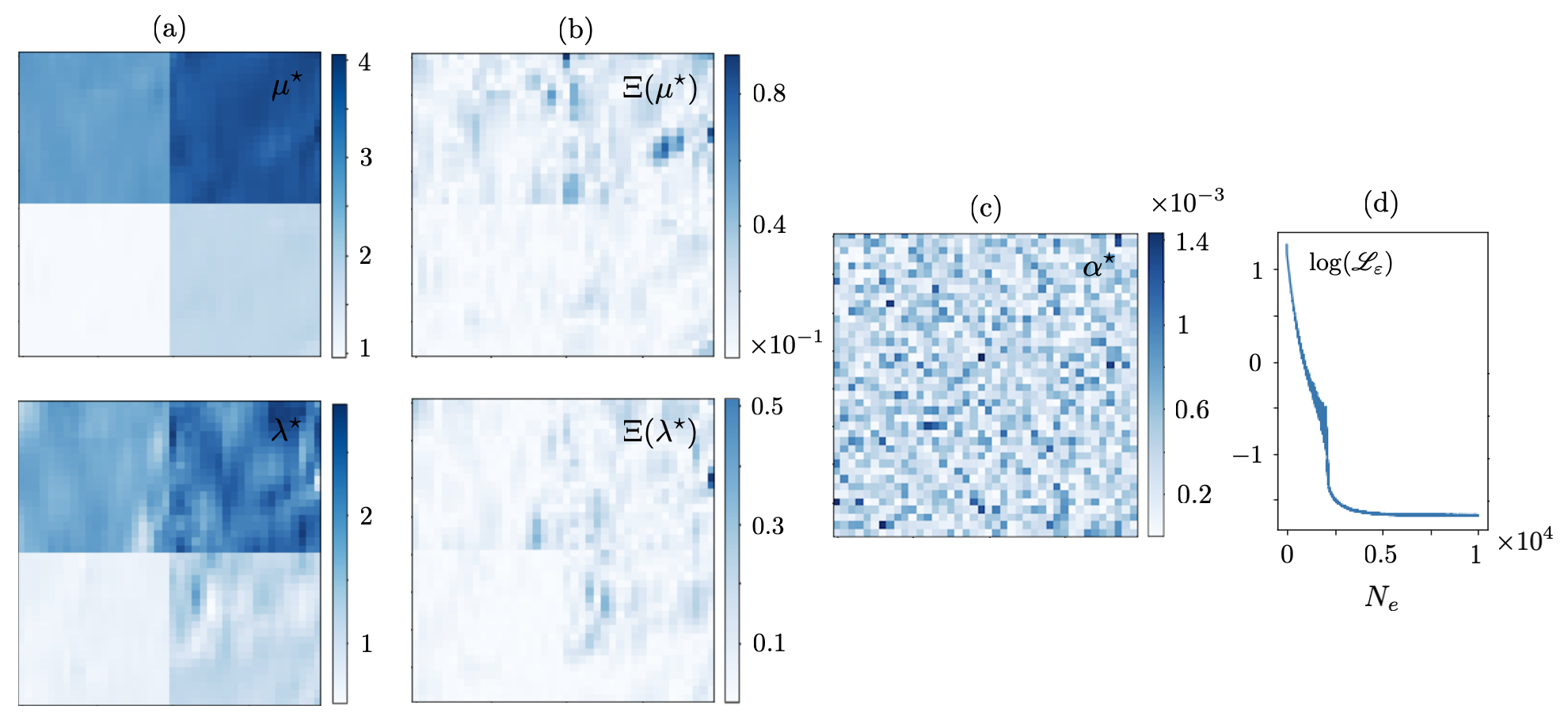} \vspace*{-2mm} 
\caption{Direct inversion of the L\'{a}me parameters in $\Uppi_2$ using five datasets perturbed with 5\% white noise:~(a)~MLP-predicted distributions $\mu^\star\nxs$ and $\lambda^\star\nxs$,~(b)~reconstruction error~\eqref{EM} with respect to the true values $\mu_j = j$ and $\lambda_j = 2j/3$, $j =  \lbrace 1, 2, 3, 4 \rbrace$,~(c) MLP-recovered regularization parameter $\alpha^\star\nxs$, and~(d) loss function $\mathcal{L}_\varepsilon$ \emph{vs}.~the number of epochs $N_e$.}
\label{DR7}
\end{figure}

To examine the impact of noise on the reconstruction, the multisource dataset used to generate Figs.~\ref{DR4} and~\ref{DR6} are perturbed with $5\%$ white noise. The subsequent direct inversions from noisy data are displayed in Figs.~\ref{DR5} and~\ref{DR7}, and the associated statistics are presented in Table~\ref{Num3}. Note that spectral differentiation as the first step in direct inversion plays a critical role in denoising the waveforms, and subsequently regularizing the reconstruction process. This may substantiate the low magnitude of MLP-recovered $\alpha^\star\nxs$ in the case of noisy data in Figs.~\ref{DR5} and~\ref{DR7}. The presence of noise, nonetheless, affects the magnitude and thus composition of terms in the Fourier representation of the processed displacement fields in space which is used for differentiation. This may in turn lead to the emergence of fluctuations in the reconstructed fields.     


\renewcommand{\arraystretch}{1.1}
\vspace*{-3mm}
 \begin{table}[!h]
\vspace*{-1mm}
\begin{center}
\caption{\small Mean $\langle \xxs \cdot \xxs \rangle_\mathcal{D}$ and standard deviation $\upsigma(\xxs \cdot \xxs |_\mathcal{D})$ of the reconstructed L\'{a}me distributions in $\mathcal{D} = \Uppi_1, \Uppi_{2_{j = 1, 2, 3, 4}}$. Here, the direct inversion is applied to noiseless data from a single source as shown in Figs.~\ref{DR15} and~\ref{DR18}.} \vspace*{2mm}
\label{Num1}
\begin{tabular}{|c|c|c|c|c|c|} \hline
\!$\mathcal{D}$\! & $\Uppi_1$\! & $\Uppi_{2_1}$\! & $\Uppi_{2_2}$\! & $\Uppi_{2_3}$\! & $\Uppi_{2_4}$\! \\ \hline\hline  
$\mu$    & $\mu_\circ\nxs = 1$  & $\mu_1\nxs = 1$  & $\mu_2\nxs = 2$ & $\mu_3\nxs = 3$ & $\mu_4\nxs = 4$ \\  
 \hline
$\exs\langle \xxs \mu^\star \rangle_\mathcal{D}$    & \!\!$0.998$\!\!  & \!\!$0.991$\!\!  & \!\!$1.983$\!\! & \!\!$2.825$\!\! & \!\!$3.835$\!\! \\ 
 \hline
$\upsigma(\mu^\star|_\mathcal{D})$   & \!\!$0.024$\!\!  & \!\!$0.083$\!\! & \!\!$0.182$\!\! & \!\!$0.441$\!\! & \!\!$0.325$\!\! \\
 \hline
$\lambda$    & $\lambda_\circ\nxs = 0.47$  & $\lambda_1\nxs = 0.67$  & $\lambda_2\nxs = 1.33$ & $\lambda_3\nxs = 2$ & $\lambda_4\nxs = 2.66$ \\  
 \hline
 $\exs\langle \xxs \lambda^\star \rangle_\mathcal{D}$    & \!\!$0.376$\!\!  & \!\!$0.615$\!\!  & \!\!$0.850$\!\! & \!\!$1.746$\!\! & \!\!$1.412$\!\! \\ 
 \hline
 $\upsigma(\lambda^\star|_\mathcal{D})$   & \!\!$0.128$\!\!  & \!\!$0.161$\!\! & \!\!$0.399$\!\! & \!\!$0.486$\!\! & \!\!$0.864$\!\! \\
 \hline
\end{tabular}
\end{center}
\vspace*{-3.5mm}
\end{table}

\renewcommand{\arraystretch}{1.1}
\vspace*{-3mm}
 \begin{table}[!h]
\vspace*{-1mm}
\begin{center}
\caption{\small Mean and standard deviation of the reconstructed L\'{a}me distributions from five distinct noiseless datasets according to Figs.~\ref{DR4} and~\ref{DR6}.} \vspace*{2mm}
\label{Num2}
\begin{tabular}{|c|c|c|c|c|c|} \hline
\!$\mathcal{D}$\! & $\Uppi_1$\! & $\Uppi_{2_1}$\! & $\Uppi_{2_2}$\! & $\Uppi_{2_3}$\! & $\Uppi_{2_4}$\! \\ \hline\hline  
$\mu$    & $1$  & $1$  & $2$ & $3$ & $4$ \\  
 \hline
$\exs\langle \xxs \mu^\star \rangle_\mathcal{D}$    & \!$1.000$\!  & \!$0.999$\!  & \!$2.003$\! & \!$2.999$\! & \!$3.999$\! \\ 
 \hline
$\upsigma(\mu^\star|_\mathcal{D})$   & \!\!$0.001$\!\!  & \!\!$0.012$\!\! & \!\!$0.011$\!\! & \!\!$0.012$\!\! & \!\!$0.016$\!\! \\
 \hline
$\lambda$    & $0.47$  & $0.67$  & $1.33$ & $2$ & $2.66$ \\  
 \hline
 $\exs\langle \xxs \lambda^\star \rangle_\mathcal{D}$    & \!\!$0.464$\!\!  & \!\!$0.660$\!\!  & \!\!$1.302$\!\! & \!\!$1.997$\!\! & \!\!$2.635$\!\! \\ 
 \hline
 $\upsigma(\lambda^\star|_\mathcal{D})$   & \!\!$0.012$\!\!  & \!\!$0.039$\!\! & \!\!$0.071$\!\! & \!\!$0.048$\!\! & \!\!$0.068$\!\! \\
 \hline
\end{tabular}
\end{center}
\vspace*{-3.5mm}
\end{table}

\renewcommand{\arraystretch}{1.1}
\vspace*{-3mm}
 \begin{table}[!h]
\vspace*{-1mm}
\begin{center}
\caption{\small Mean and standard deviation of the reconstructed L\'{a}me distributions from noisy data according to Figs.~\ref{DR5} and~\ref{DR7}.}
\label{Num3}
\begin{tabular}{|c|c|c|c|c|c|} \hline
\!$\mathcal{D}$\! & $\Uppi_1$\! & $\Uppi_{2_1}$\! & $\Uppi_{2_2}$\! & $\Uppi_{2_3}$\! & $\Uppi_{2_4}$\! \\ \hline\hline  
$\mu$    & $1$  & $1$  & $2$ & $3$ & $4$ \\   
 \hline
$\exs\langle \xxs \mu^\star \rangle_\mathcal{D}$    & $1.001$  & $1.002$  & $2.005$ & $2.996$ & $3.996$ \\ 
 \hline
$\upsigma(\mu^\star|_\mathcal{D})$   & \!\!$0.005$\!\!  & \!\!$0.016$\!\! & \!\!$0.035$\!\! & \!\!$0.054$\!\! & \!\!$0.088$\!\! \\
 \hline
$\lambda$    & $0.47$  & $0.67$  & $1.33$ & $2$ & $2.66$ \\  
 \hline
 $\exs\langle \xxs \lambda^\star \rangle_\mathcal{D}$    & \!\!$0.462$\!\!  & \!\!$0.650$\!\!  & \!\!$1.263$\!\! & \!\!$2.006$\!\! & \!\!$2.654$\!\! \\ 
 \hline
 $\upsigma(\lambda^\star|_\mathcal{D})$   & \!\!$0.042$\!\!  & \!\!$0.051$\!\! & \!\!$0.225$\!\! & \!\!$0.182$\!\! & \!\!$0.300$\!\! \\
 \hline
\end{tabular}
\end{center}
\vspace*{-3.5mm}
\end{table}

\subsection{Physics-informed neural networks}\label{PINN-synt}

The learning process of Section~\ref{PINN} is performed as follows:~(a)~the MLP network ${\boldsymbol{u}^\upalpha}^\star = \mathcal{N}_{\boldsymbol{u}^\upalpha}(\bxi,\omega,\bx \exs | \exs \gamma,\boldsymbol{\vartheta}^\star)$ endowed with the positive-definite parameters $\gamma$ and $\boldsymbol{\vartheta}^\star = (\mu^\star, \lambda^\star)$ is constructed such that the input $\bx$ labels the source location and the auxiliary weight ${\gamma}$ is a nonadaptive scaler,~(b)~$\mu^\star$ and $\lambda^\star$ may be specified as scaler or distributed parameters of the network according to Fig.~\ref{DIPINN4}~(i), and~(c)~${\boldsymbol{u}^\upalpha}^\star$ is trained by minimizing $\mathcal{L}_\varpi$ in~\eqref{LPINNs} via the ADAM optimizer using the synthetic waveforms of Section~\ref{PS-synt}. Reconstructions are performed on the same set of collocation points sampling $S^\text{obs} \!\times\nxs \Omega \nxs \times\! \mathcal{T}$ as in Section~\ref{DI-synt}. Accordingly, the input to ${\boldsymbol{u}^\upalpha}^\star$ is of size $N_\xi \!\times\! N_\omega \!\times\! N_\tau = 1600 \!\times\! 1 \!\times\! N_s$, while its output is of dimension $(1600 \!\times\! 1 \!\times\! N_s)^2$ modeling the displacement field along $\xi_1$ and $\xi_2$ in the sampling region. Similar to Section~\ref{DI-synt}, each epoch makes use of the full dataset for training and the learning rate is $0.005$. The PyTorch implementation of PINNs in this section is accomplished by building upon the available codes on the Github repository~\cite{Roxis2022}.   

The MLP network ${\boldsymbol{u}^1}^\star = {\boldsymbol{u}^1}^\star (\bxi,\omega,\bx \exs | \exs \gamma,\boldsymbol{\vartheta}^\star)$ with three hidden layers of respectively $20$, $40$, and  $20$ neurons is employed to map the displacement field ${\boldsymbol{u}^1}$ (in $\Uppi_1$) associated with a single point source of frequency $\omega = 3.91$ at $\bx = \bx_1 \in S^{\text{inc}}$. The L\'{a}me constants are defined as the unknown \emph{scaler} parameters of the network i.e.,~$\boldsymbol{\vartheta}^\star = \lbrace \mu^\star, \lambda^\star \rbrace$, and the Lagrange multiplier $\gamma$ is specified per the following argument. Within the dimensional framework of this section and with reference to~\eqref{map-synt}, observe that on setting $\gamma = \frac{1}{\rho \omega^2}$ (i.e., $\gamma = 0.065$), both (the PDE residue and data misfit) components of the loss function $\mathcal{L}_{\varpi}$ in~\ref{LPINNs} emerge as some form of balance in terms of the displacement field. This may naturally facilitate maintaining of the same scale for the loss terms during training, and thus, simplifying the learning process by dispensing with the need to tune an additional parameter $\gamma$. Keep in mind that the input to ${\boldsymbol{u}^1}^\star$ is of size $1600 \!\times\! 1 \!\times\! 1$, while its output is of dimension $(1600 \!\times\! 1 \!\times\! 1)^2$. In this setting, the training objective is two-fold:~(a) construction of a surrogate map for ${\boldsymbol{u}^1}$, and~(b) identification of $\mu^\star$ and $\lambda^\star$. 

\begin{figure}[!h]
\center\includegraphics[width=0.95\linewidth]{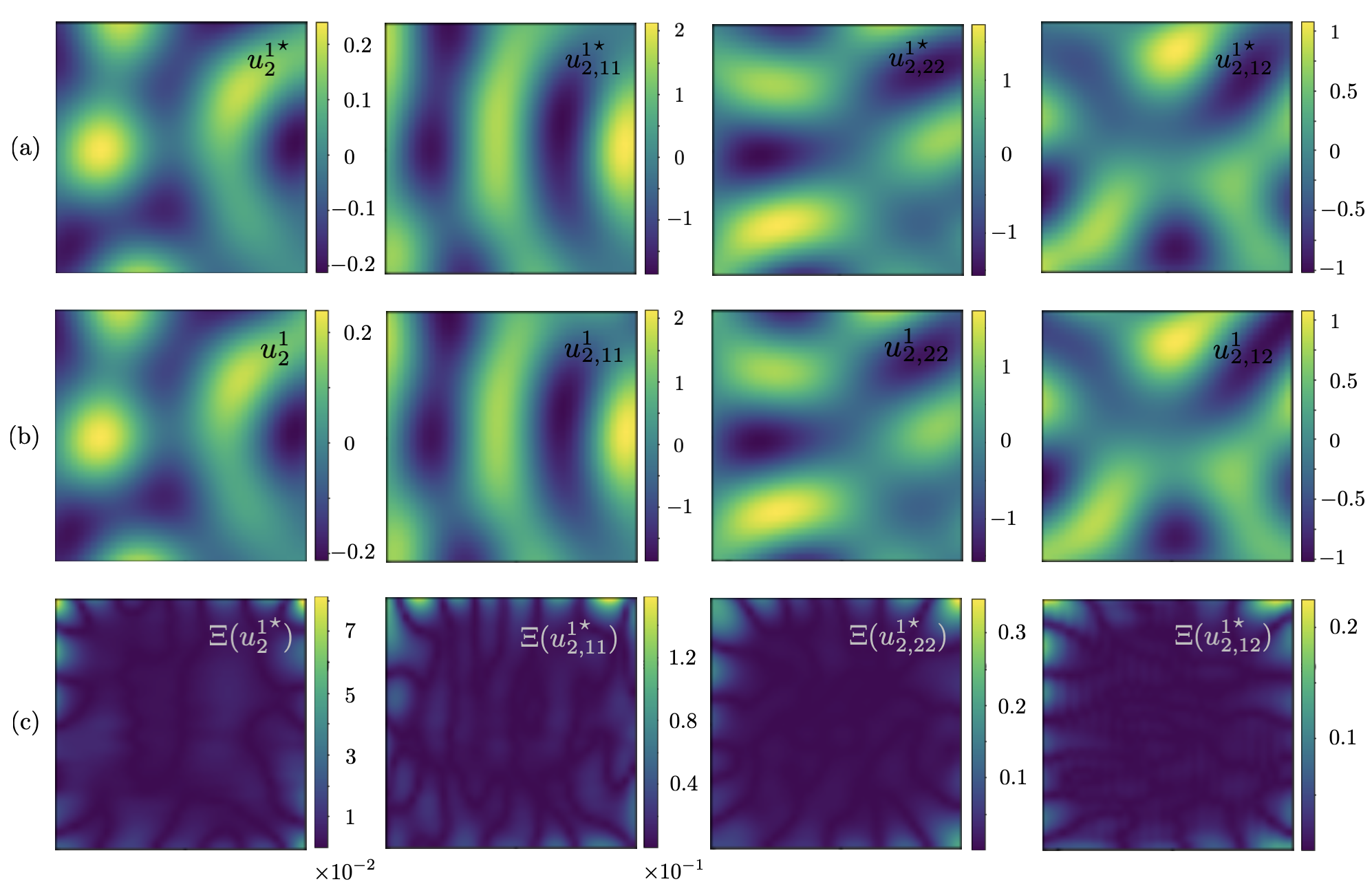} \vspace*{-4mm} 
\caption{PINN \emph{vs}.~FEM maps of vertical displacement and its derivatives in $\Uppi_1$:~(a) MLP estimates, from noiseless data, for $\lbrace {u_2^1}^\star\!, {u^1}^\star_{\!\!\!\hspace{-0.1mm} 2,11}, {u^1}^\star_{\!\!\!\hspace{-0.1mm} 2,22}, {u^1}^\star_{\!\!\!\hspace{-0.1mm} 2,12} \rbrace$ wherein the derivatives ${u^1}^\star_{\!\!\!\hspace{-0.1mm} 2,ij}$, $i,j = 1,2$, are obtained by automatic differentiation, (b) FEM displacement solution and its spectral derivatives for $\lbrace {u}_2^1, {u}^1_{2,11}, {u}^1_{2,22}, {u}^1_{2,12} \rbrace$, and~(c) normal misfit~\ref{EM} between (a) and (b). }
\label{PINNuy}
\vspace*{-1mm} 
\end{figure}

Fig.~\ref{PINNuy} showcases (i) the accuracy of PINN estimates based on noiseless data in terms of the vertical component of displacement field ${u}^1_2$ in $\Uppi_1$, and (ii) the performance of automatic differentiation~\cite{pasz2017} in capturing the field derivatives in terms of components that appear in the governing PDE~\ref{map-synt} i.e., ${u}^1_{2,ij} = {\partial^2 {u}^1_2}/(\partial \xi_i \partial \xi_j)$, $i,j = 1,2$. The comparative analysis in (ii) is against the spectral derivates of FEM fields according to Section~\ref{DIT}. It is worth noting that similar to Fourier-based differentiation, the most pronounced errors in automatic differentiation occur in the near-boundary region i.e., the support of one-sided derivatives. It is observed that the magnitude of such discrepancies may be reduced remarkably by increasing the number of epochs. Nonetheless, the loci of notable errors remain at the vicinity of specimen's external boundary or internal discontinuities such as cracks or material interfaces. Fig.~\ref{PINNuy} is complemented with the reconstruction       

\begin{figure}[!h]
\center\includegraphics[width=0.91\linewidth]{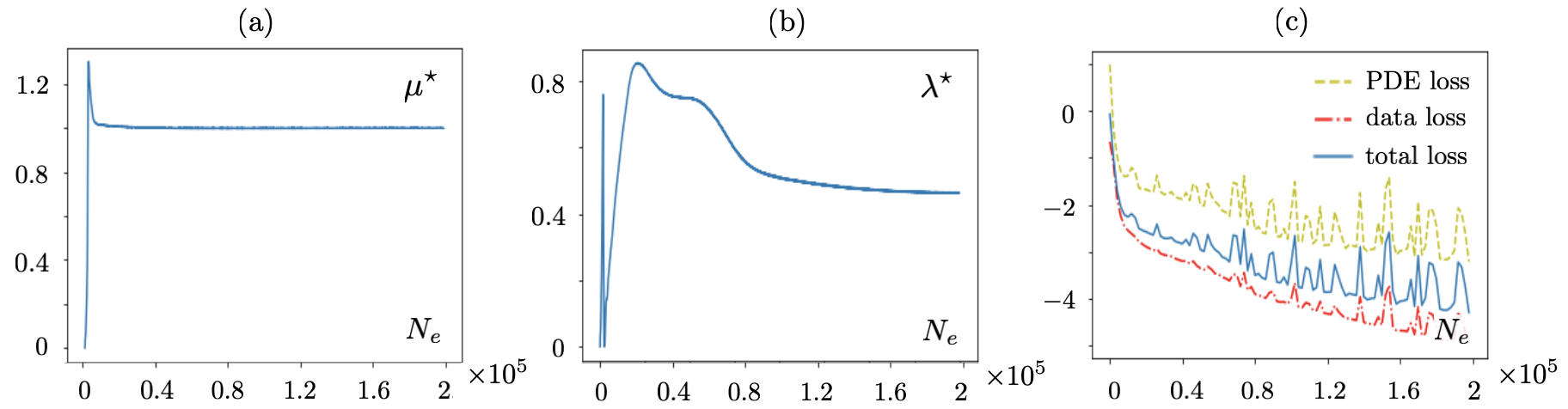} \vspace*{-2mm} 
\caption{PINN reconstruction of L\'{a}me constants in the homogeneous plate $\Uppi_1$ from noiseless data:~(a) $\mu^\star$ \emph{vs}. number of epochs $N_e$,~(b) $\lambda^\star$ \emph{vs}. $N_e$, and (c) total loss $\mathcal{L}_{\varpi}$ and its components (the PDE residue and data misfit) \emph{vs.} $N_e$ in log scale.}
\label{PINNHNL}
\end{figure} 
\begin{figure}[!h]
\center\includegraphics[width=0.915\linewidth]{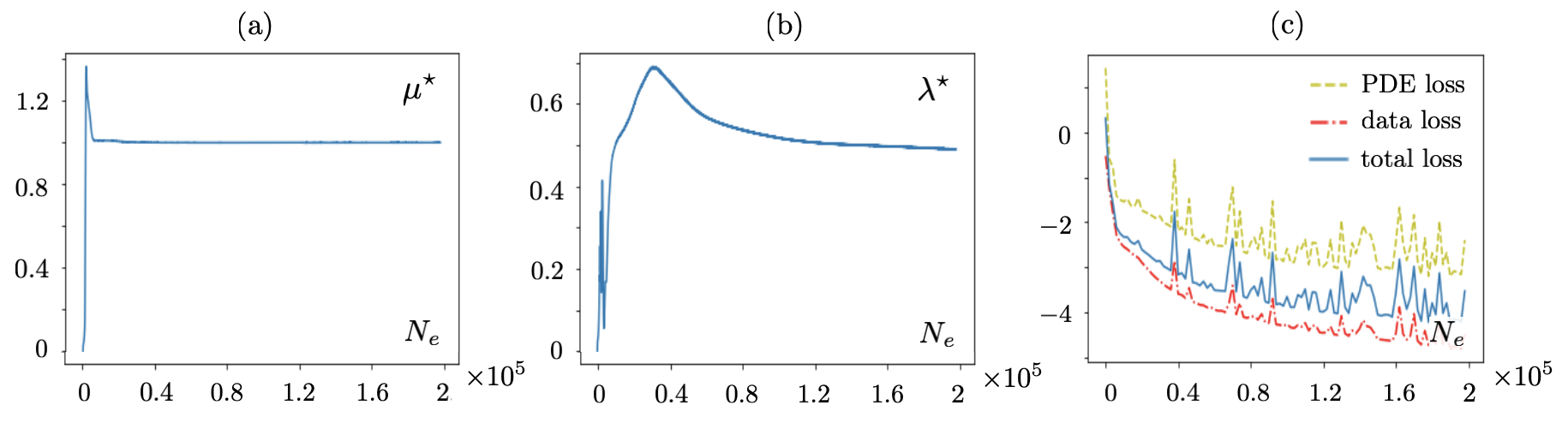} \vspace*{-4mm} 
\caption{PINN reconstruction of L\'{a}me constants in $\Uppi_1$ from noisy data:~(a) $\mu^\star$ \emph{vs}. number of epochs $N_e$,~(b) $\lambda^\star$ \emph{vs}. $N_e$, and (c) total loss $\mathcal{L}_{\varpi}$ and its components (the PDE residue and data misfit) \emph{vs.} $N_e$ in log scale.}
\label{DRP2}
\vspace*{-4mm}
\end{figure} 

\noindent results of Fig.~\ref{PINNHNL} indicating $(\mu^\star, \lambda^\star) = (1.000, 0.486)$ for the homogenous specimen $\Uppi_1$ with the true L\'{a}me constants $(\mu_\circ, \lambda_\circ) = (1, 0.47)$. The impact of noise on training is examined by perturbing the noiseless data related to Fig.~\ref{PINNuy} with $5\%$ white noise, which led to $(\mu^\star, \lambda^\star) = (0.999, 0.510)$ as shown in Fig.~\ref{DRP2}.                 

\begin{figure}[!h]
\center\includegraphics[width=0.95\linewidth]{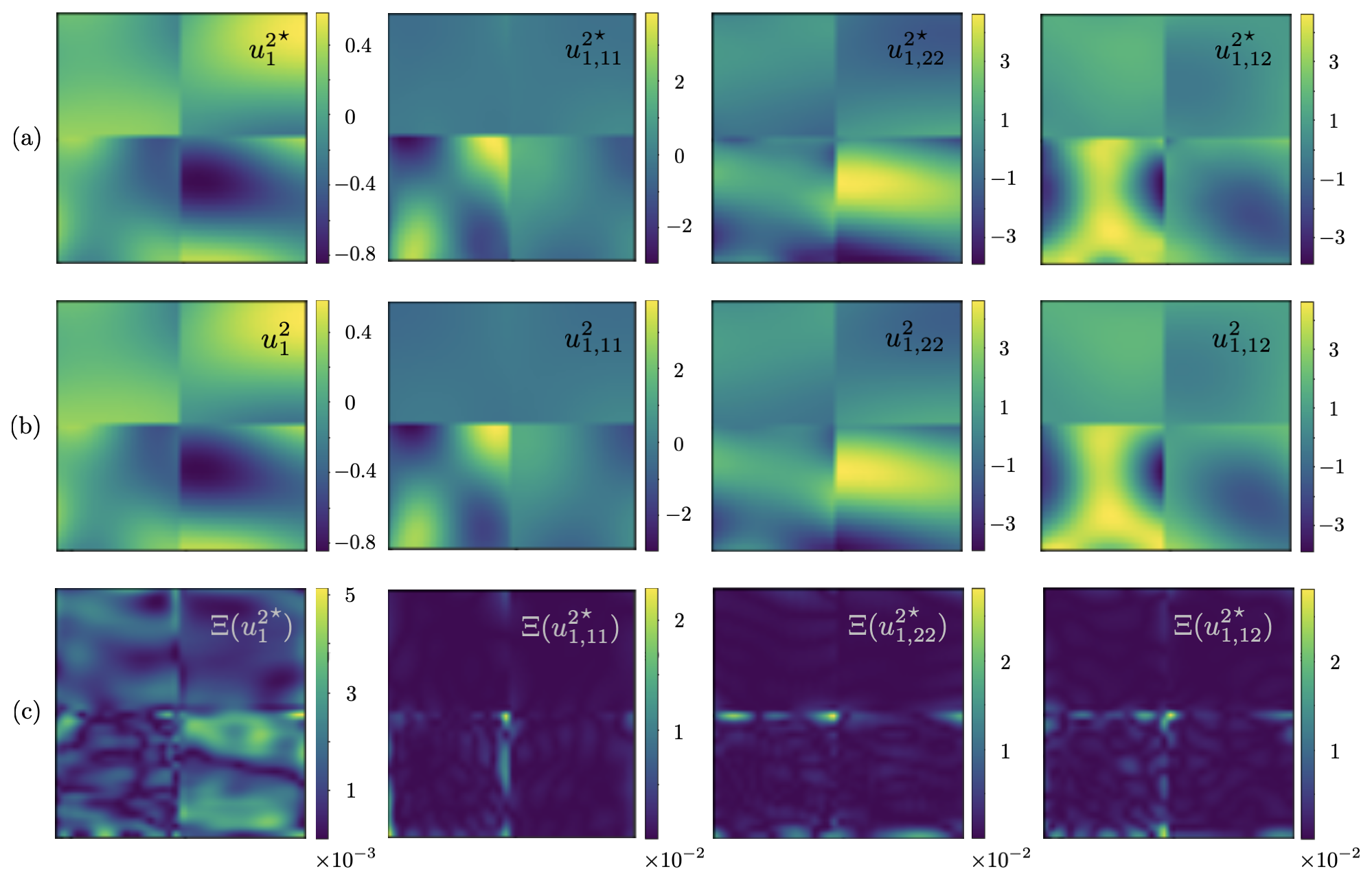} \vspace*{-3mm} 
\caption{PINN \emph{vs}.~FEM maps of horizontal displacement and its derivatives in $\Uppi_2$:~(a) PINN estimates, from noiseless data, for $\lbrace {u_1^2}^\star\!, {u^2}^\star_{\!\!\!\hspace{-0.1mm} 1,11}, {u^2}^\star_{\!\!\!\hspace{-0.1mm} 1,22}, {u^2}^\star_{\!\!\!\hspace{-0.1mm} 1,12} \rbrace$ wherein the derivatives ${u^2}^\star_{\!\!\!\hspace{-0.1mm} 1,ij}$, $i,j = 1,2$, are obtained by automatic differentiation, (b) FEM displacement solution and its spectral derivatives for $\lbrace {u}_1^2, {u}^2_{1,11}, {u}^2_{1,22}, {u}^2_{1,12} \rbrace$, and~(c) normal misfit~\ref{EM} between (a) and (b). }
\label{PINNux}
\vspace*{-2mm}
\end{figure}

Next, the PINN ${\boldsymbol{u}^2}^\star = {\boldsymbol{u}^2}^\star (\bxi,\omega,\bx \exs | \exs \boldsymbol{\vartheta}^\star)$ with three hidden layers of respectively $120$, $120$, and $80$ neurons  is created to reconstruct (i) displacement field ${\boldsymbol{u}^2}$ in the heterogeneous specimen $\Uppi_2$, and (ii) distribution of the L\'{a}me parameters over the observation surface. In this vein, synthetic waveform data associated with five point sources $\lbrace \bx_i \rbrace \in S^{\text{inc}}$, $i = 1, 2, \ldots, 5$ at $\omega = 3.91$ is used for training. Here, $\boldsymbol{\vartheta}^\star$ is the network's unknown \emph{distributed} parameter, of dimension $(40 \!\times\! 40)^2$, and the \emph{nonadaptive scaler} weight $\gamma = 0.065$ in light of the sample's uniform density $\rho = 1$. In this setting, the input to ${\boldsymbol{u}^2}^\star$ is of size $1600 \!\times\! 1 \!\times\! 5$, while its output is of dimension $(1600 \!\times\! 1 \!\times\! 5)^2$. Fig.~\ref{PINNux} provides a comparative analysis between the FEM and PINN maps of horizontal displacement ${u}^1_2$ in $\Uppi_2$ and its spatial derivatives computed by spectral and automatic differentiation respectively.

\begin{figure}[!h]
\center\includegraphics[width=0.84\linewidth]{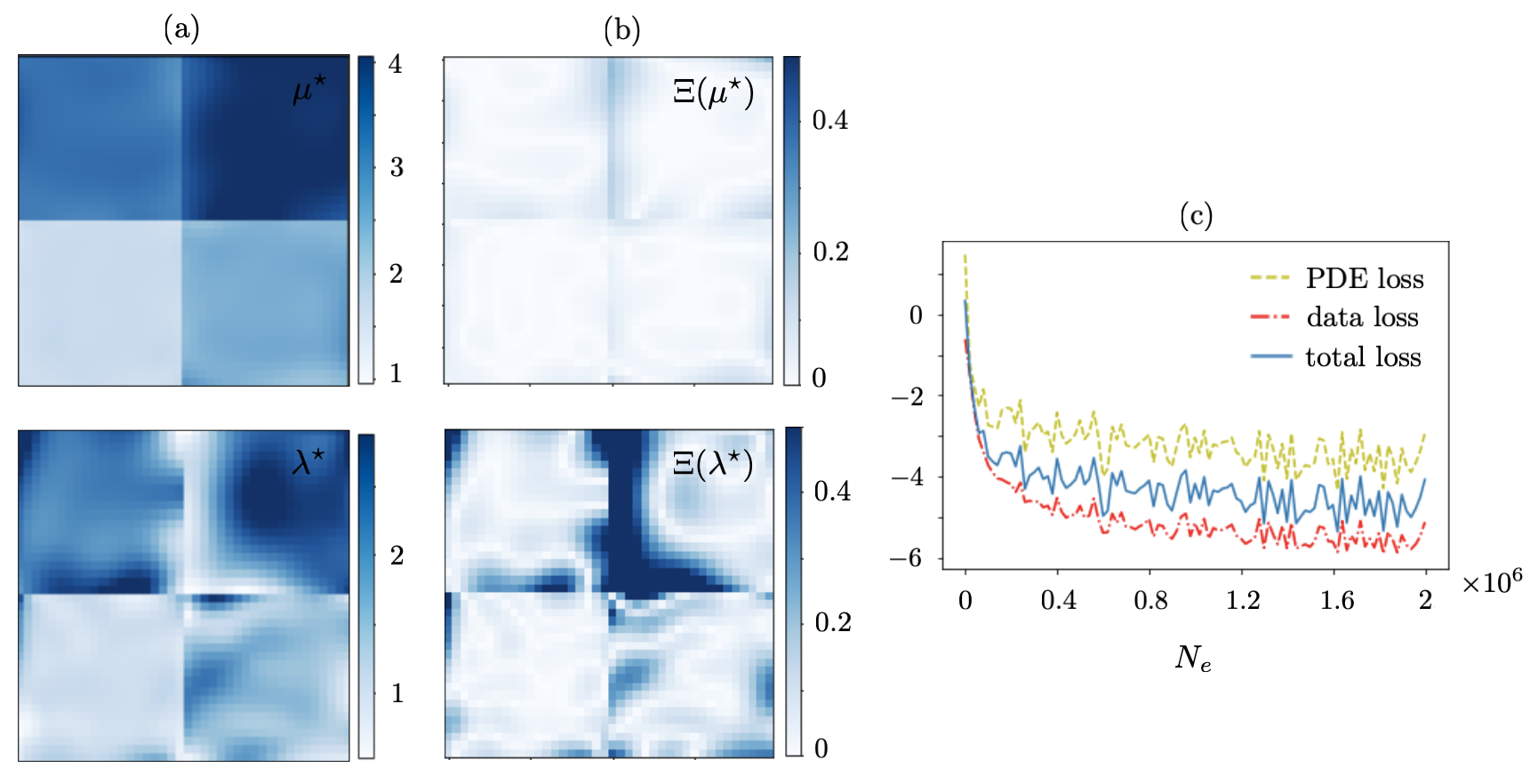} \vspace*{-3mm} 
\caption{PINN reconstruction of L\'{a}me parameters in $\Uppi_2$ using five noiseless datasets:~(a)~PINN-predicted distributions $\mu^\star\nxs$ and $\lambda^\star\nxs$,~(b)~reconstruction error~\eqref{EM} with respect to the true values $\mu_j = j$ and $\lambda_j = 2j/3$, $j =  \lbrace 1, 2, 3, 4 \rbrace$,~(c) total loss $\mathcal{L}_{\varpi}$ and its components (the PDE residue and data misfit) \emph{vs.} $N_e$ in log scale.}
\label{PINNHNL2}
\vspace*{-2mm}
\end{figure} 

\renewcommand{\arraystretch}{1.1}
\vspace*{-3mm}
 \begin{table}[!h]
\vspace*{-1mm}
\begin{center}
\caption{\small Mean and standard deviation of the PINN-reconstructed L\'{a}me distributions from five distinct noiseless datasets according to Fig.~\ref{PINNHNL2}.} \vspace*{2mm}
\label{Num4}
\begin{tabular}{|c|c|c|c|c|} \hline
\!$\mathcal{D}$\! & $\Uppi_{2_1}$\! & $\Uppi_{2_2}$\! & $\Uppi_{2_3}$\! & $\Uppi_{2_4}$\! \\ \hline\hline  
$\exs\langle \xxs \mu^\star \rangle_\mathcal{D}$   & $0.975$  & $1.973$ & $2.941$ &. $3.918$ \\ 
 \hline
$\upsigma(\mu^\star|_\mathcal{D})$  & \!\!$0.054$\!\! & \!\!$0.123$\!\! & \!\!$0.135$\!\! & \!\!$0.226$\!\! \\
 \hline
$\exs\langle \xxs \lambda^\star \rangle_\mathcal{D}$   & \!\!$0.686$\!\!  & \!\!$1.250$\!\! & \!\!$2.045$\!\! & \!\!$2.065$\!\! \\ 
 \hline
 $\upsigma(\lambda^\star|_\mathcal{D})$   & \!\!$0.247$\!\! & \!\!$0.400$\!\! & \!\!$0.520$\!\! & \!\!$0.857$\!\! \\
 \hline
\end{tabular}
\end{center}
\vspace*{-3.5mm}
\end{table}

The PINN-reconstructed distribution of PDE parameters is illustrated in Fig.~\ref{PINNHNL2} whose statistics is detailed in Table~\ref{Num4}. It is worth mentioning that the learning process is repeated for a suit of weights $\gamma = \lbrace 0.01, 0.025, 0.1, 0.25, 0.5, 1.5, 2, 5, 10, 15 \rbrace$. In all cases, the results are either quite similar or worse than that of Figs.~\ref{PINNux} and~\ref{PINNHNL2}.

\section{Laboratory implementation}\label{Le}

This section examines the performance of direct inversion and PINNs for full-field ultrasonic characterization in a laboratory setting. In what follows, experimental data are processed prior to inversion as per Section~\ref{DIT} which summarizes the detailed procedure  in~\cite{pour2018(2)}. To verify the inversion results, quantities of interest are also reconstructed through dispersion analysis, separately, from a set of auxiliary experiments. 

\subsection{Test set-up}\label{PS-exp}

\begin{figure}[!bp]
\vspace*{-4mm}
\center\includegraphics[width=0.9\linewidth]{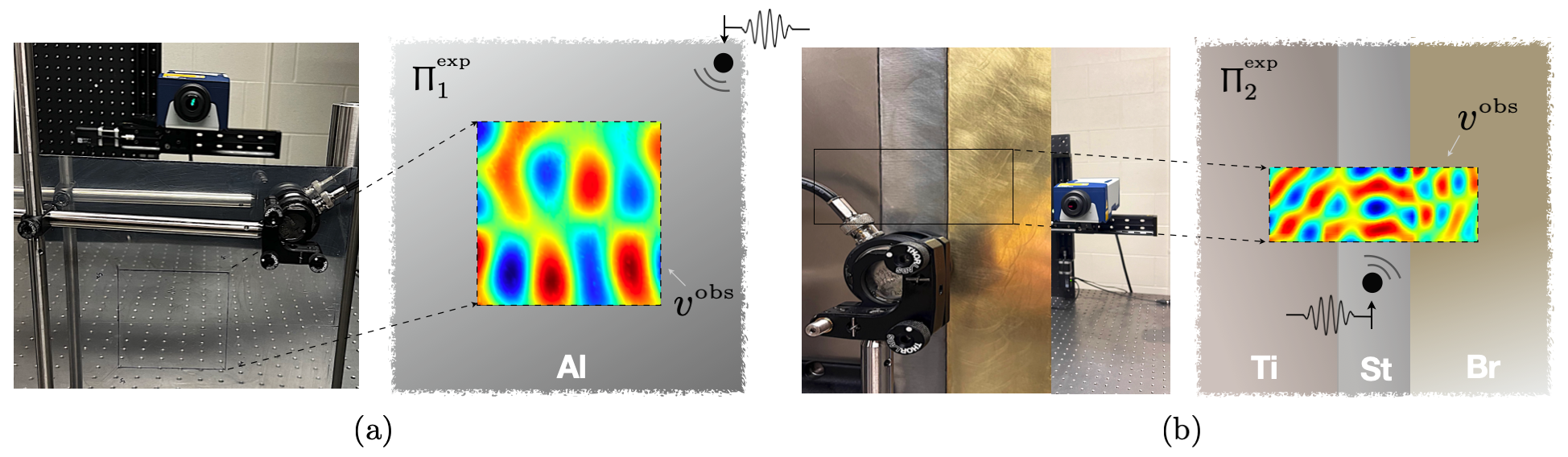} \vspace*{-3mm} 
\caption{Test set-ups for ultrasonic full-field characterization:~(a)~an Al plate $\Uppi_1^{^\text{exp}}\!$ is subject to antiplane shear waves at $165$ kHz by a piezoelectric transducer; the out-of-plane particle velocity field is then captured by a laser Doppler vibrometer scanning on a robot over the observation surface, and~(b) a Ti-St-Br plate $\Uppi_2^{^\text{exp}}\!$ undergoes a similar test at $80$ kHz and $300$ kHz.}
\label{SE3}
\vspace*{-4mm}
\end{figure} 

Experiments are performed on two (homogeneous and heterogeneous) specimens: $\Uppi_1^{^\text{exp}}\!\nxs$ which is a $27$~cm $\!\nxs\times 27$ cm $\!\nxs\times 1.5$ mm sheet of T6 6061 aluminum, and $\Uppi_2^{^\text{exp}}\!$ composed of~(a)~$5$ cm $\!\nxs\times\xxs 27$ cm $\!\nxs\times 1.5$ mm sheet of Grade 2 titanium,~(b)~$2.5$ cm $\!\nxs\times\xxs 27$ cm $\!\nxs\times 1.5$ mm sheet of 4130 steel, and~(c)~$5$ cm $\!\nxs\times\xxs 27$ cm $\!\nxs\times 1.5$ mm sheet of 260-H02 brass, connected via metal epoxy. For future reference, the density $\uprho_{\mu}$, Young's modulus $\textsf{E}_{\mu}$, and Poisson's ratio $\nu_{\mu}$ for $\mu = \lbrace \text{Al, Ti, St, Br} \rbrace$ are listed in Table~\ref{prop} as per the manufacturer.      

Ultrasonic experiments on both samples are performed in a similar setting in terms of the sensing configuration and illuminating wavelet. In both cases, the specimen is excited by an \emph{antiplane} shear wave from a designated source location $S^{\text{inc}}$, shown in Fig.~\ref{SE3}, by a $0.5$ MHz p-wave piezoceramic transducer (V101RB by Olympus Inc.). The incident signal is a five-cycle burst of the form 
\begin{equation}\label{wavelet}
H({\sf f_c t}) \, H(5\!-\!{\sf f_c t}) \, \sin\big(0.2 \pi {\sf f_c t}\big) \, \sin\big(2 \pi {\sf f_c t}\big), 
\end{equation}
where $H$ denotes the Heaviside step function, and the center frequency ${\sf f_c}\!$ is set at $165$ kHz (\emph{resp}.~$\lbrace 80, 300\rbrace$ kHz) in $\Uppi_1^{^\text{exp}}\!$ (\emph{resp}.~$\Uppi_2^{^\text{exp}}\nxs$). The induced wave motion is measured in terms of the particle velocity ${\sf{v}^\upbeta}$, $\upbeta = 1,2$, on the scan grids $\mathcal{G}_\upbeta$ sampling $S^{\text{obs}}$ where $S^{\text{obs}} \cap S^{\text{inc}} = S^{\text{obs}} \cap \partial\Uppi_\upbeta^{^\text{exp}}\! = \emptyset$. A laser Doppler vibrometer (LDV) which is mounted on a 2D robotic translation frame (for scanning) is deployed for measurements. The VibroFlex Xtra VFX-I-120 LDV system by Polytec Inc.~is capable of capturing particle velocity within the frequency range $\sim\nxs\text{DC}-24$ MHz along the laser beam which in this study is normal to the specimen's surface. 

The scanning grid $\mathcal{G}_1 \subset \Uppi_1^{^\text{exp}}\!$ is identified by a $2$ cm $\!\nxs\times\xxs 2$ cm square sampled by $100 \!\times\! 100$ uniformly spaced measurement points. This amounts to a spatial resolution of $0.2$ mm in both spatial directions. In parallel, $\mathcal{G}_2 \subset \Uppi_2^{^\text{exp}}\!$ is a $2.5$ cm $\!\nxs\times\xxs 7.5$ cm rectangle positioned according to Fig.~\ref{SE3}~(b) and sampled by a uniform grid of $180 \!\times\! 60$ scan points associated with the spatial resolution of $0.42$ mm. At every scan point, the data acquisition is conducted for a time period of $400$ $\mu$s at the sampling rate of $250$ MHz. To minimize the impact of optical and mechanical noise in the system, the measurements are averaged over an ensemble of 80 realizations at each scan point. Bear in mind that both the direct inversion and PINNs deploy the spectra of normalized \emph{velocity} fields $v^\text{obs}$ for data inversion. Such distributions of out-of-plane particle velocity at $165$ kHz (\emph{resp}.~$80$ kHz) in $\Uppi_1^{^\text{exp}}\!$ (\emph{resp}.~$\Uppi_2^{^\text{exp}}\nxs$) is displayed in Fig.~\ref{SE3}. 

It should be mentioned that in the above experiments, the magnitude of measured signals in terms of displacement is of $O(\text{nm})$ so that it may be appropriate to assume a linear regime of propagation. The nature of antiplane wave motion is dispersive nonetheless. Therefore, to determine the relevant length scales in each component, the associated dispersion curves are obtained as in Fig.~\ref{Disp} via a set of complementary experiments described in Section~\ref{VAV-exp}. Accordingly, for excitations of center frequency $\lbrace {\sf f_{c_1}}, {\sf f_{c_2}}, {\sf f_{c_3}} \rbrace  = \lbrace 165, 80, 300 \rbrace$ kHz, the affiliated phase velocity $\textsf{c}_{\mu}$ and wavelength $\uplambda_{\mu}$ for $\mu = \lbrace \text{Al, Ti, St, Br} \rbrace$ is identified in Table~\ref{prop3}.  

\subsection{Dimensional framework}\label{DP-exp}

On recalling Section~\ref{DP}, let $\ell_r \colon\!\!\!= \uplambda_{\text{Al}} = 0.01$ m, $\mu_r \colon\!\!\!= \textsf{E}_{\text{Al}} = 68.9$ GPA, and $\rho_r  \colon\!\!\!= \uprho_{\text{Al}} = 2700$ kg/m$^3$ be the reference scales for length, stress, and mass density, respectively. In this setting, the following maps take the physical quantities to their dimensionless values
\beq\lb{Dmap}
\begin{aligned}
& (\uprho_{\mu}, \textsf{E}_{\mu}, \nu_{\mu}) ~\rightarrow~ (\rho_{\mu}, E_{\mu}, \nu_{\mu}) := \big(\dfrac{1}{\rho_r}\uprho_{\mu}, \dfrac{1}{\mu_r}\textsf{E}_{\mu}, \nu_{\mu} \big), & \!\!\!\!\!\! \mu = \lbrace \text{Al, Ti, St, Br} \rbrace, \\*[0.5mm]
& ({\sf f_{c_\iota}}, \uplambda_{\mu}, \textsf{c}_{\mu}) ~\rightarrow~ ({ f_{c_\iota}}, \lambda_{\mu}, {c}_{\mu}) := \big(\ell_r\sqrt{\frac{\rho_r}{\mu_r}}\exs {\sf f_{c_\iota}}, \dfrac{1}{\ell_r}\uplambda_{\mu}, \sqrt{\frac{\rho_r}{\mu_r}} \exs \textsf{c}_{\mu} \big), &\!\!\!\!\!\! \iota = 1,2,3, \\*[0.25mm]
& ({\sf h}, {\sf f}, {\sf{v}^\upbeta}) ~\rightarrow~ (h, f, v^\upbeta) := \big(\dfrac{1}{\ell_r} {\sf h}, \ell_r\sqrt{\frac{\rho_r}{\mu_r}}\exs {\sf f}, \sqrt{\frac{\rho_r}{\mu_r}}\exs {\sf{v}^\upbeta} \big), & \!\!\!\!\!\! \upbeta = 1,2,
\end{aligned}
\vspace*{-1mm} 
\eeq
where ${\sf h} = 1.5$ mm and ${\sf f}$ respectively indicate the specimen's thickness and cyclic frequency of wave motion. Table~\ref{prop} (\emph{resp}.~Table~\ref{prop3}) details the normal values for the first (\emph{resp}.~second) of~\eqref{Dmap}. The normal thickness and center frequencies are as follows,   
\beq\lb{Dmap2}
\lbrace {f_{c_1}}, {f_{c_2}}, {f_{c_3}} \rbrace  = \lbrace 0.33, 0.16, 0.59 \rbrace, \quad h = 0.15.
\eeq

\vspace*{-4.5mm}
\renewcommand{\arraystretch}{1.1}
 \begin{table}[!h]
\vspace*{-1mm}   
 \begin{center}
 \caption{\small Properties of the aluminum, titanium, steel and brass sheets as per the manufacturer. Here, $\chi_{\mu}\colon \!\!\!= {E_\mu}/{\rho_\mu}$.} \vspace*{1mm}
\label{prop}
 \begin{tabular}{c c c c c c}
\hline
\multirow{2}{*}{\vspace*{-5mm}{\small{\bf physical}}} & \!\small{$\mu$}\!                    & \!\small{Al}\! & \!\small{Ti}\! & \!\small{St}\! & \!\small{Br}\!   \\\cline{2-6}
                                                           & \!\small{$\textsf{E}_{\mu}$ \footnotesize{[GPA]}}\! & \!\small{68.9}\!                  & \!\small{105}\!                 & \!\small{199.95}\!                   & \!\small{110}\!                      \\\cline{2-6}
\multirow{2}{*}{\vspace*{4.5mm}{\small{\bf quantity}}}     & \!\small{$\uprho_{\mu}$ \!\footnotesize{[kg/m$^3$]}}\! & \!\small{2700}\!                  & \!\small{4510}\!                 & \!\small{7850}\!                   & \!\small{8530}\!                      \\\cline{2-6}
                                                           & \!\small{$\nu_{\mu}$}\! & \!\small{0.33}\!                  & \!\small{0.34}\!                 & \!\small{0.29}\!                   & \!\small{0.31}\!                      \\\hline          
\multirow{3}{*}{\vspace*{4mm}{\small{\bf normal}}}   & \!\small{$E_{\mu}$}\! & \!\small{1}\!                  & \!\small{1.52}\!                 & \!\small{2.90}\!                   & \!\small{1.60}\!                      \\\cline{2-6}
\multirow{3}{*}{\vspace*{5mm}{\small{\bf value}}}                                                           & \!\small{$\rho_{\mu}$}\! & \!\small{1}\!                  & \!\small{1.67}\!                 & \!\small{2.91}\!                   & \!\small{3.16}\!                      \\\cline{2-6}     
                                                                                  & \!\small{$\chi_{\mu}$}\! & \!\small{1}\!                  & \!\small{0.91}\!                 & \!\small{1}\!                   & \!\small{0.51}\!                      \\\hline                                                                                                                                                                      
\end{tabular}
\end{center}
\vspace*{-4.5mm}
\end{table}

\vspace*{-5.5mm}
\begin{table}[!htb]
    \caption{\small Phase velocity $\textsf{c}_{\mu}$ and wavelength $\uplambda_{\mu}$ in $\mu = \lbrace \text{Al, Ti, St, Br} \rbrace$ at $\lbrace {\sf f_{c_1}}, {\sf f_{c_2}}, {\sf f_{c_3}} \rbrace  = \lbrace 165, 80, 300 \rbrace$ kHz as per Fig.~\ref{Disp}, and their normalized counterparts according to~\eqref{Dmap}.}\vspace*{2.5mm}
\label{prop3}  
    \begin{minipage}{.55\linewidth}
   {\small{\bf physical quantity}}\vspace*{1mm}
      \centering
 \begin{tabular}{c c c c c}  
\hline
 \!\small{$\mu$}\!        & \!\small{Al}\!             & \!\small{Ti}\! & \!\small{St}\! & \!\small{Br}\!   \\\hline
 \!\small{$\uplambda_{\mu}({\sf f_{c_1\!}})$ \footnotesize{[cm]}}\! & \!\small{$1$}\!   & \!\small{$-$}\!                  & \!\small{$-$}\!                 & \!\small{$-$}\!                                       \\\hline
 \!\small{$\textsf{c}_{\mu}({\sf f_{c_1\!}})$ \!\footnotesize{[m/s]}}\! & \!\small{$1610.4$}\!  & \!\small{$-$}\!                  & \!\small{$-$}\!                 & \!\small{$-$}\!                                       \\\hline   
\!\small{$\uplambda_{\mu}({\sf f_{c_2\!}})$ \footnotesize{[cm]}}\! & \!\small{$-$}\!   & \!\small{$1.4$}\!                  & \!\small{$1.4$}\!                 & \!\small{$1.17$}\!                                       \\\hline
 \!\small{$\textsf{c}_{\mu}({\sf f_{c_2\!}})$ \!\footnotesize{[m/s]}}\!  & \!\small{$-$}\!     & \!\small{$1140$}\!                  & \!\small{$1126$}\!                 & \!\small{$936$}\!                                   \\\hline   
 \!\small{$\uplambda_{\mu}({\sf f_{c_3\!}})$ \footnotesize{[cm]}}\! & \!\small{$-$}\!   & \!\small{$0.65$}\!                  & \!\small{$0.64$}\!                 & \!\small{$0.5$}\!                                       \\\hline
 \!\small{$\textsf{c}_{\mu}({\sf f_{c_3\!}})$ \!\footnotesize{[m/s]}}\! &\!\small{$-$}\!    & \!\small{$1960.8$}\!                  & \!\small{$1929$}\!                 & \!\small{$1501.6$}\!                                      \\\hline                                                                                                                                                                         
\end{tabular}
    \end{minipage}%
    \begin{minipage}{.38\linewidth}
      \centering    
{\small{\bf normal value}}\vspace*{1mm}      
 \begin{tabular}{c c c c c}
\hline
 \!\small{$\mu$}\!        & \!\small{Al}\!             & \!\small{Ti}\! & \!\small{St}\! & \!\small{Br}\!   \\\hline
 \!\small{$\lambda_{\mu}({ f_{c_1\!}})$}\!  & \!\small{$1$}\!  & \!\small{$-$}\!                  & \!\small{$-$}\!                 & \!\small{$-$}\!                                        \\\hline
 \!\small{${c}_{\mu}({f_{c_1\!}})$}\! & \!\small{$0.32$}\!   & \!\small{$-$}\!                  & \!\small{$-$}\!                 & \!\small{$-$}\!                                       \\\hline   
\!\small{$\lambda_{\mu}({f_{c_2\!}})$}\! & \!\small{$-$}\! & \!\small{$1.4$}\!                  & \!\small{$1.4$}\!                 & \!\small{$1.17$}\!                                         \\\hline
 \!\small{${c}_{\mu}({f_{c_2\!}})$}\! & \!\small{$-$}\!    & \!\small{$0.23$}\!                  & \!\small{$0.22$}\!                 & \!\small{$0.19$}\!                                      \\\hline   
 \!\small{$\lambda_{\mu}({f_{c_3\!}})$}\! & \!\small{$-$}\!   & \!\small{$0.65$}\!                  & \!\small{$0.64$}\!                 & \!\small{$0.5$}\!                                       \\\hline
 \!\small{${c}_{\mu}({f_{c_3\!}})$}\! & \!\small{$-$}\! & \!\small{$0.39$}\!                  & \!\small{$0.38$}\!                 & \!\small{$0.3$}\!                                        \\\hline                                                                                                                                                                         
\end{tabular}
    \end{minipage} 
\end{table}

\subsection{Governing equation}\label{GE-exp}

In light of~\eqref{Dmap2} and Table~\ref{prop3}, observe that in all tests the wavelength-to-thickness ratio $\frac{\lambda_{\mu}}{h} \in [3.33 \,\,\, 9.33]$, $\mu = \lbrace \text{Al, Ti, St, Br} \rbrace$. Therefore, one may invoke the equation governing flexural waves in thin plates~\cite{graf2012} to approximate the physics of measured data. In this framework,~\eqref{uf} may be recast as 
\beq\label{map-exp}
\begin{aligned}
&\!\!  \Lambda ~=~\Lambda_\upbeta ~\colon\!\!\!=~ \frac{ \chi_\upbeta h^3}{12(1-\nu_{\nxs\upbeta}^2)}\nabla^4   \,-\,  h (2 \pi {f})^2 ,  & \chi_\upbeta \,\colon \!\!\!=\, \frac{E_\upbeta}{\rho_\upbeta}, \, \upbeta \,=\, 1,2,  \\*[0.45mm] 
& \!\! \hat{\text{{u}}} ~=~ v^\upbeta(\bxi,f;\btau), \quad \boldsymbol{\vartheta} ~=~ \chi_\upbeta(\bxi,f),  &  \bxi \in S^\text{obs}\!, \, \btau \in S^{\text{inc}}, \, f \in [0.8 \,\,\, 1.2]f_{c_\iota}, \, \iota = 1,2,3,
\end{aligned}
\eeq
where $\rho_\upbeta, E_\upbeta, \nu_{\nxs\upbeta}$ respectively denote the normal density, Young's modulus, and Poisson's ratio in $\Uppi_\upbeta^{^\text{exp}}\!\nxs$, $\upbeta \,=\, 1,2$, and $\btau$ indicates the source location. Note that $\nu_{\nxs\upbeta} \sim 0.32$ according to Table~\ref{prop} and $\Lambda$, related to $1-\nu_{\nxs\upbeta}^2$, shows little sensitivity to small variations in the Poisson's ratio. Thus, in what follows, $\nu_{\nxs\upbeta}$ is treated as a \emph{known} parameter. Provided $v^\upbeta(\bxi,f;\btau)$, the objective is to reconstruct $\chi_\upbeta(\bxi,f)$.

\subsection{Direct inversion}\label{DI-exp}

Following the reconstruction procedure of Section~\ref{DI-synt}, the distribution of $\chi_\upbeta$ in $\mathcal{G}_\upbeta$, $\upbeta = 1,2$, is obtained at specific frequencies. In this vein, the positive-definite MLP networks $\chi^\star_\upbeta = \mathcal{N}_{\chi_\upbeta}(\bxi,\omega)$ and $\alpha^\star = \mathcal{N}_{\alpha}(\bxi,\omega)$ comprised of three hidden layers of respectively $20$, $40$, and $20$ neurons are constructed according to~Fig.~\ref{DIDM7}. In all MLP trainings of this section, each epoch makes use of the full dataset and the learning rate is $0.005$. 

When $\upbeta = 1$, the inversion is conducted at $f_1 = 0.336$. $S^\text{inc}$ is sampled at one point i.e., the piezoelectric transducer remains fixed during the test on Al plate, and thus, $N_\tau =1$, while a concentric $60\!\times\!60$ subset of collocation points sampling $S^{\text{obs}}$ is deployed for training. In this setting, the input to $\chi^\star_1$ and $\alpha^\star$ is of size $N_\xi N_\tau \nxs\times\nxs N_\omega = 3600 \nxs\times\nxs 1$, and their real-valued outputs are of the same size. The results are shown in Fig.~\ref{DR8}. When $\upbeta = 2$, the direct inversion is conducted at $f_2 = 0.17$ and $f_3 = 0.61$. For the low-frequency reconstruction, $S^\text{inc}$ is sampled at one point, while a $40\!\times\!120$ subset of scan points in $\mathcal{G}_2$ is used for training so that the input/output size for $\chi^\star_2$ and $\alpha^\star$ is $4600 \!\times\! 1$. The recovered fields and associated normal error are provided in Fig.~\ref{DR20}. Table~\ref{Num5} enlists the true values as well as mean and standard deviation of the reconstructed distributions $\chi^\star_\upbeta$ in $\Uppi_\upbeta^{^\text{exp}}\!\nxs$, $\upbeta = 1,2$, according to Figs.~\ref{DR8} and~\ref{DR20}. For the high-frequency reconstruction, when $\upbeta = 2$, $S^\text{inc}$ is sampled at three points i.e., experiments are performed for three distinct positions of the piezoelectric transducer, while the same subset of scan points is used for training. In this case, the input to $\chi^\star_{2}$ and $\alpha^\star$ is $13800 \!\times\! 1$, while their output is of dimension $4600 \!\times\! 1$. The high-frequency reconstruction results are illustrated in Fig.~\ref{DR21}, and the affiliated means and standard deviations are provided in Table~\ref{Num6}. It should be mentioned that the computed normal errors in Figs.~\ref{DR8},~\ref{DR20}, and~\ref{DR21} are with respect to the verified values of Section~\ref{VAV-exp}. Note that the recovered $\alpha^\star$s from laboratory test data are much smoother than the ones reconstructed from synthetic data in Section~\ref{DI-synt}. This could be attributed to the scaler nature of~\eqref{map-exp} with a single unknown parameter -- as opposed to the vector equations governing the in-plane wave motion with two unknown parameters. More specifically, here, $\alpha^\star$ controls the weights and biases of a single network $\chi^\star_\upbeta$, while in Section~\ref{DI-synt}, $\alpha^\star$ simultaneously controls the parameters of two separate networks $\mu^\star$ and $\lambda^\star$. A comparative analysis of Figs.~\ref{DR20} and~\ref{DR21} reveals that~(a) enriching the waveform data by increasing the number of sources remarkably decrease the reconstruction error,~(b)~the regularization parameter $\alpha$ in~\eqref{Lels} is truly distributed in nature as the magnitude of the recovered $\alpha^\star$ in brass is ten times greater than that of titanium and steel which is due to the difference in the level of noise in measurements related to distinct material surfaces, and (c) the recovered field $\chi^\star_{2}$ -- which according to~\eqref{map-exp} is a material property ${E_2}/{\rho_2}$, demonstrates a significant dependence to the reconstruction frequency. The latter calls for proper verification of the results which is the subject of Section~\ref{VAV-exp}.   

\subsubsection{Verification}\label{VAV-exp}

To shine some light on the nature discrepancies between the low- and high- frequency reconstructions in
 
\begin{figure}[!bp]
 \vspace*{-4mm} 
\center\includegraphics[width=0.94\linewidth]{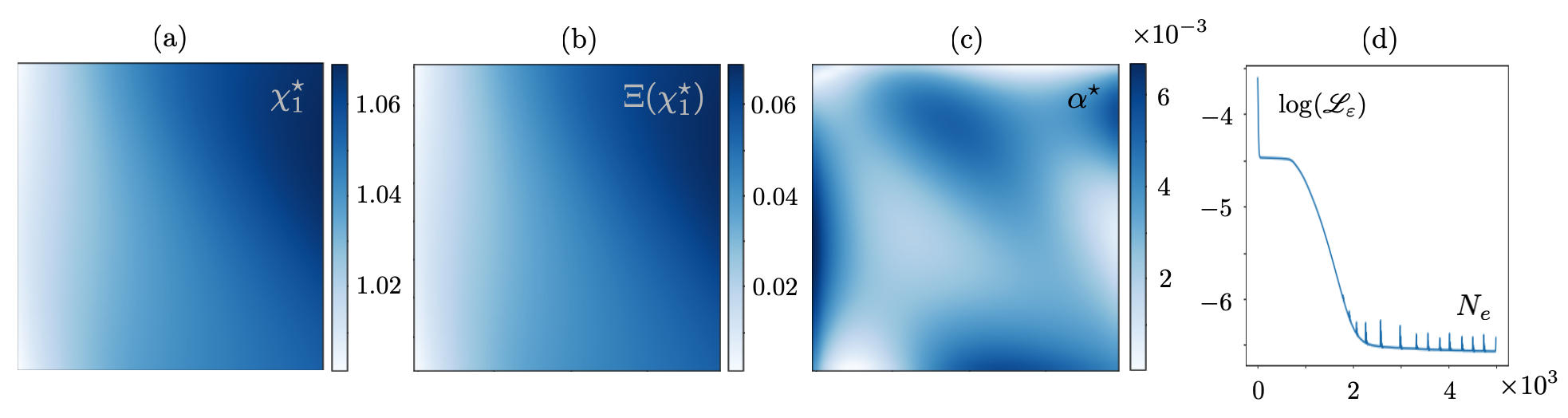} \vspace*{-4mm} 
\caption{Direct inversion of the PDE parameter $\chi_1$ in $\Uppi_1^{^\text{exp}}\!\nxs$ using test data from a single source at frequency $f_1 = 0.336$:~(a)~MLP-predicted distribution $\chi_1(\bxi,f_1)$ in $\bxi \in \mathcal{G}_1$,~(b)~reconstruction error~\eqref{EM} with respect to the true value $\chi_1 = \chi_{\text{Al}} =1$,~(c) MLP-recovered distribution of the regularization parameter $\alpha^\star\nxs$, and~(d) loss function $\mathcal{L}_\varepsilon$ \emph{vs}.~the number of epochs $N_e$ in log scale.}
\label{DR8}
\end{figure} 

\pagebreak

\begin{figure}[!h]
\center\includegraphics[width=0.952\linewidth]{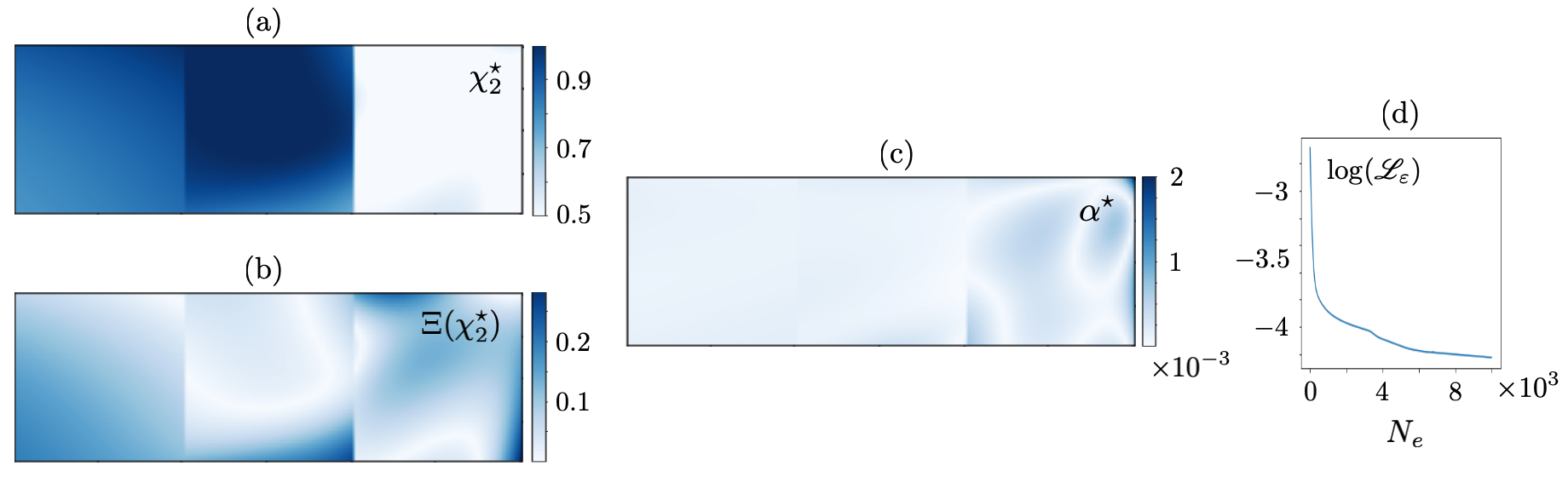} \vspace*{-3mm} 
\caption{Direct inversion of the PDE parameter $\chi_2$ in $\Uppi_2^{^\text{exp}}\!\nxs$ using test data from a single source at frequency $f_2 = 0.17$:~(a)~MLP-predicted distribution $\chi_2(\bxi,f_2)$ in $\bxi \in \mathcal{G}_2$,~(b)~reconstruction error~\eqref{EM} with respect to the true value $\chi_2 \in \lbrace \chi_{\text{Ti}}, \chi_{\text{St}}, \chi_{\text{Br}} \rbrace = \lbrace 0.91, 1, 0.51 \rbrace$ as per Table~\ref{prop},~(c) MLP-recovered distribution of the regularization parameter $\alpha^\star\nxs$, and~(d) loss function $\mathcal{L}_\varepsilon$ \emph{vs}.~the number of epochs $N_e$ in log scale.}
\vspace*{-2mm}
\label{DR20}
\end{figure} 

\renewcommand{\arraystretch}{1.1}
\vspace*{-3mm}
 \begin{table}[!h]
\vspace*{-1mm}
\begin{center}
\caption{\small Mean and standard deviation of the reconstructed distributions in Figs.~\ref{DR8} and~\ref{DR20} via the direct inversion of single-source test data.} \vspace*{2mm}
\label{Num5}
\begin{tabular}{|c|c|c|c|c|} \hline
\!$\upbeta$\! & $1$\! & $2_{\text{Ti}}$\! & $2_{\text{St}}$\! & $2_{\text{Br}}$\!  \\ \hline\hline  
$\chi_\upbeta$    & $1$  & $0.91$  & $1$ & $0.51$  \\   
 \hline
$\exs\langle \xxs \chi_\upbeta^\star \rangle_{\Uppi_\upbeta^{^\text{exp}}\!\nxs}$    & $1.041$  & $0.872$  & $0.978$ & $0.443$ \\ 
 \hline
$\upsigma(\chi_\upbeta^\star|_{\Uppi_\upbeta^{^\text{exp}}\!\nxs})$   & \!\!$0.017$\!\!  & \!\!$0.044$\!\! & \!\!$0.060$\!\! & \!\!$0.052$\!\!  \\
 \hline
\end{tabular}
\end{center}
\vspace*{-3.5mm}
\end{table}

\begin{figure}[!h]
\center\includegraphics[width=0.943\linewidth]{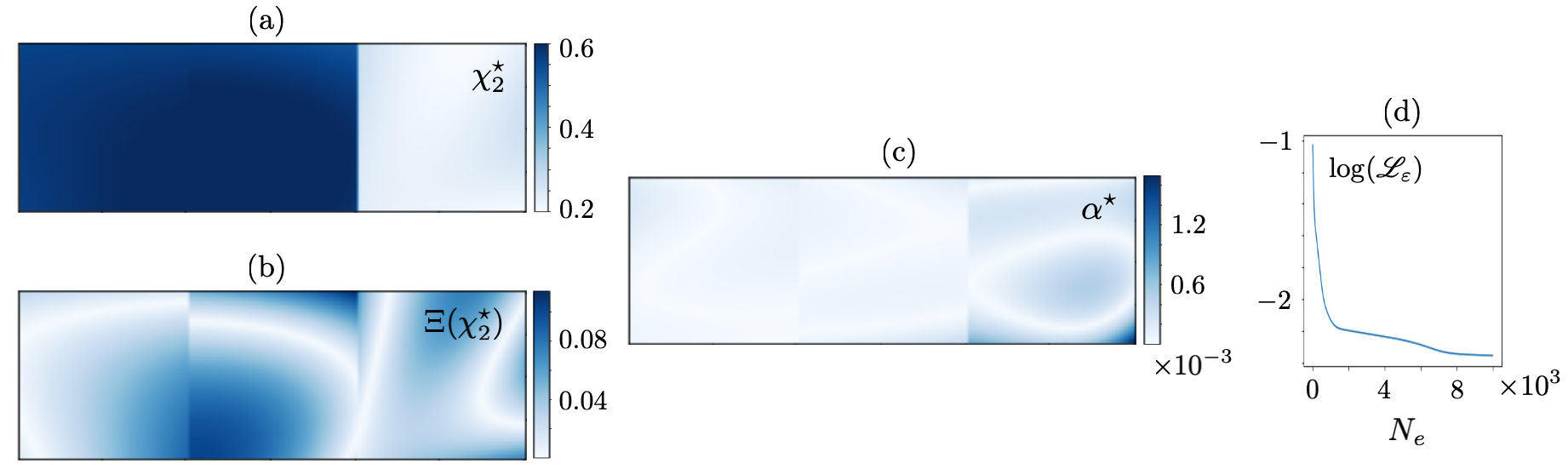} \vspace*{-3mm} 
\caption{Direct inversion of the PDE parameter $\chi_2$ in $\Uppi_2^{^\text{exp}}\!\nxs$ using test data from three source locations at frequency $f_3 = 0.61$:~(a)~MLP-predicted distribution $\chi_2(\bxi,f_3)$ in $\bxi \in \mathcal{G}_2$,~(b)~reconstruction error~\eqref{EM} with respect to the related estimates $\lbrace 0.57, 0.59, 0.24 \rbrace$ as per Fig.~\ref{DR22},~(c) MLP-recovered distribution of the regularization parameter $\alpha^\star\nxs$, and~(d) loss function $\mathcal{L}_\varepsilon$ \emph{vs}.~the number of epochs $N_e$ in log scale.}
\label{DR21}
\end{figure} 

\renewcommand{\arraystretch}{1.1}
\vspace*{-3mm}
 \begin{table}[!h]
\vspace*{-1mm}
\begin{center}
\caption{\small Mean and standard deviation of the reconstructed distributions in Fig.~\ref{DR21} via the direct inversion applied to high-frequency test data from three distinct sources.} \vspace*{2mm}
\label{Num6}
\begin{tabular}{|c|c|c|c|} \hline
\!$\upbeta$\!  & $2_{\text{Ti}}$\! & $2_{\text{St}}$\! & $2_{\text{Br}}$\!  \\ \hline\hline  
$\chi'_\upbeta$    & $0.57$  & $0.59$ & $0.24$  \\   
 \hline
$\exs\langle \xxs \chi_\upbeta^\star \rangle_{\Uppi_\upbeta^{^\text{exp}}\!\nxs}$   & \!\!$0.585$\!\!  & \!\!$0.606$\!\! & \!\!$0.227$\!\! \\ 
 \hline
$\upsigma(\chi_\upbeta^\star|_{\Uppi_\upbeta^{^\text{exp}}\!\nxs})$   & \!\!$0.015$\!\! & \!\!$0.029$\!\! & \!\!$0.016$\!\!  \\
 \hline
\end{tabular}
\end{center}
\end{table}

\begin{figure}[!tp]
\center\includegraphics[width=0.97\linewidth]{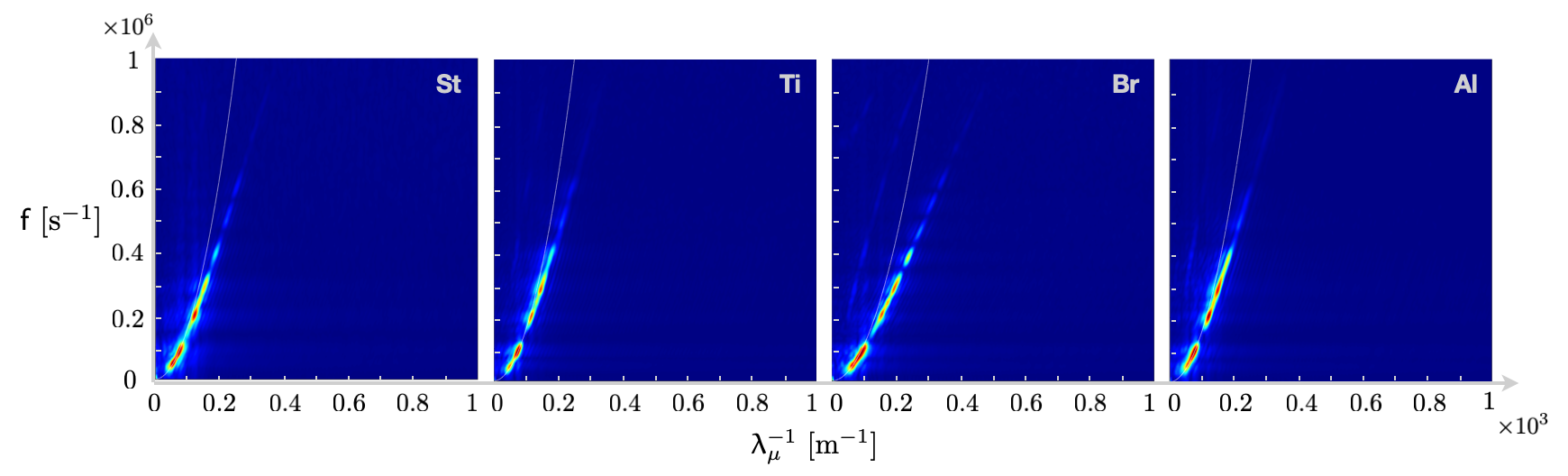} \vspace*{-4mm} 
\caption{Experimental \emph{vs}.~theoretical dispersion curves ${\sf f}(\uplambda_\mu^{-1})$ for $\mu = \lbrace \text{Al, Ti, St, Br} \rbrace$. Analytical curves (solid lines) are computed from~\eqref{tdisp} using the pertinent properties in Table~\ref{prop}.}\label{Disp}
\end{figure} 

\begin{figure}[!h]
\center\includegraphics[width=0.85\linewidth]{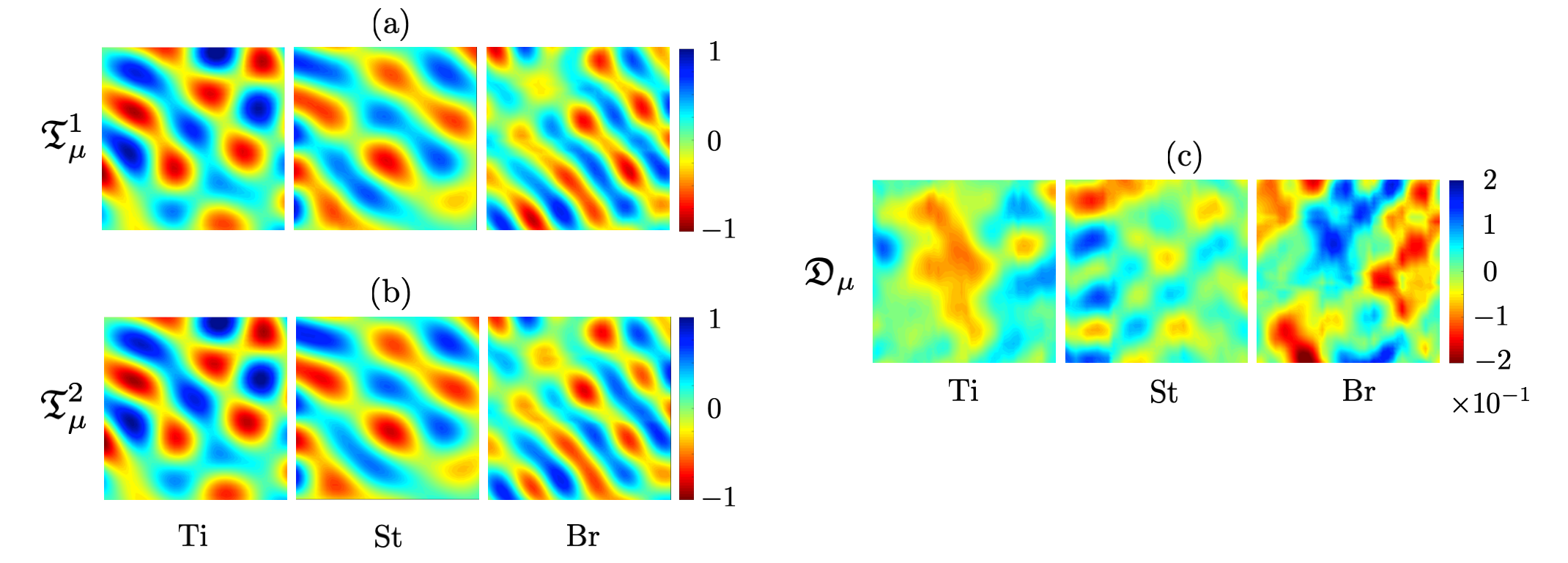} \vspace*{-3mm} 
\caption{Discrepancy in the balance law~\eqref{map-exp} at $f_3 = 0.61$:~(a)~elastic force field $\mathfrak{T}_\mu^1$, $\mu = \lbrace \text{Ti, St, Br} \rbrace$, according to~\eqref{map-exp2} with adjusted coefficients $\lbrace \chi_{\text{Ti}}, \chi_{\text{St}}, \chi_{\text{Br}} \rbrace = \lbrace 0.57, 0.59, 0.24 \rbrace$,~(b) the inertia field $\mathfrak{T}_\mu^2$, and~(c) normal discrepancy $\mathfrak{D}_\mu$.}
\label{DR22}
\vspace*{-2.5mm} 
\end{figure} 

\noindent Figs.~\ref{DR20} and~\ref{DR21}, a set of secondary tests are performed to obtain the dispersion curve for each component of the test setup. For this purpose, antiplane shear waves of form~\eqref{wavelet} are induced at ${\sf f_{c_j}}\! = 50 j$ kHz, $j = 1,2, \ldots, 7$, in $60$ cm $\!\times\!$ $60$ cm cuts of aluminum, titanium, steel, and brass sheets used in the primary tests of Fig.~\ref{SE3}. In each experiment, the piezoelectric transducer is placed in the middle of specimen (far from the external boundary), and the out-of-plane wave motion is captured in the immediate vicinity of the transducer along a straight line of length $8$ cm sampled at $400$ scan points. The Fourier-transformed signals in time-space furnish the dispersion relations of Fig.~\ref{Disp}. In parallel, the theoretical dispersion curves affiliated with~\eqref{map-exp} are computed according to
\beq\lb{tdisp}
{\sf f}\! ~=~ 2\pi(\uplambda_\mu)^{-2} \sqrt{\frac{\chi_\mu {\sf h}^2}{12 (1-\nu_\mu^2)}}, \quad  \chi_\mu ~=~ \frac{{\sf E}_\mu}{\uprho_\mu}, \quad \mu = \lbrace \text{Al, Ti, St, Br} \rbrace,
\eeq
using the values of Table~\ref{prop} for $\chi_\mu$ and $\nu_\mu$ and ${\sf h} = 1.5$mm. A comparison between the experimental and theoretical dispersion curves ${\sf f}(\uplambda_\mu^{-1})$ in Fig.~\ref{Disp} verifies the theory and the values of Table~\ref{prop} for $\chi_\mu$ in the low-frequency regime of wave motion. This is also in agreement with the direct inversion results of Figs.~\ref{DR8} and~\ref{DR20}. Moreover, Fig.~\ref{Disp} suggests that at approximately ${\sf f}_\mu = \lbrace{170, 200, 120, 110} \rbrace$ kHz for $\mu = \lbrace \text{Al, Ti, St, Br} \rbrace$ the governing PDE~\eqref{map-exp} with physical coefficients fails to predict the experimental results which may provide an insight regarding the high-frequency reconstruction results in Fig.~\ref{DR21}. Further investigation of the balance law~\eqref{map-exp}, as illustrated in Fig.~\ref{DR22}, shows that the test data at $312$ kHz satisfy -- with less than $10-20\%$ discrepancy depending on the material -- a PDE of form~\eqref{map-exp} with modified coefficients. More specifically, Fig.~\ref{DR22} demonstrates the achievable balance between the elastic force distribution $\mathfrak{T}_\mu^1$ and inertia field $\mathfrak{T}_\mu^2$ in~\eqref{map-exp} by directly adjusting the PDE parameter $\chi'_2$ to minimize the discrepancy $\mathfrak{D}_\mu$ according to  

\beq\label{map-exp2}
\begin{aligned}
&\!\!  \mathfrak{T}_\mu^1 ~\colon\!\!\!=~ \frac{ \chi'_2 h^3}{12(1-\nu_{\nxs 2}^2)}\nabla^4 v^2, \quad   \mathfrak{T}_\mu^2 ~\colon\!\!\!=~  h (2 \pi {f})^2 v^2, \quad \mathfrak{D}_\mu ~\colon\!\!\!=~ \frac{|\mathfrak{T}_\mu^1 - \mathfrak{T}_\mu^2|}{\max|\mathfrak{T}_\mu^2|}.
\end{aligned}
\eeq
With reference to Table~\ref{Num6}, the recovered coefficients $\chi'_2$ at $f = f_3 = 0.61$ verify the direct inversion results of Fig.~\ref{DR21}. This implies that the direct inversion (or PINNs) may lead to non-physical reconstructions in order to attain the best fit for the data to the ``perceived"" underlying physics. Thus, it is imperative to establish the range of validity of the prescribed physical principles in data-driven modeling. Here, the physics of the system at $f_3$ is in transition, yet close enough to the leading-order approximation~\eqref{map-exp} that the discrepancy is less than $20\%$. It is unclear, however, if this equation with non-physical coefficients may be used as a predictive tool. It would be interesting to further investigate the results through the prism of higher-order continuum theories and a set of independent experiments for validation which could be the subject of a future study.

\subsection{Physics-informed neural networks}\label{PINN-exp}

Following Section~\ref{PINN-synt}, PINNs are built and trained using experimental test data of Section~\ref{DI-exp}. The MLP network ${v^1}^\star = {v^1}^\star (\bxi,f,\bx \exs | \exs \gamma,\chi_1^\star)$ with six hidden layers of respectively $40$, $40$, $120$, $80$, $40$, and $40$ neurons is constructed to map the out-of-plane velocity field ${v^1}$ (in $\Uppi_1^{^\text{exp}}\nxs$) related to a single transducer location $\bx_1$ and frequency $f_1 = 0.336$. The PDE parameter $\chi_1$ is defined as the unknown scaler parameter of the network, and following the argument of Section~\ref{PINN-synt}, the Lagrange multiplier $\gamma$ is specified as a nonadaptive scaler weight of magnitude $\frac{1}{h (2\pi f_1)^2} = 1.5$. The input/output dimension for ${v^1}^\star$ is $N_\xi \!\times\! N_\omega \!\times\! N_\tau = 3600 \!\times\! 1 \!\times\! 1$, and each epoch makes use of the full dataset for training and the learning rate is $0.005$. Keep in mind that the objective here is to~(a) construct a surrogate map for ${v^1}$, and~(b) identify $\chi_1^\star$. 

\begin{figure}[!tp]
\center\includegraphics[width=0.95\linewidth]{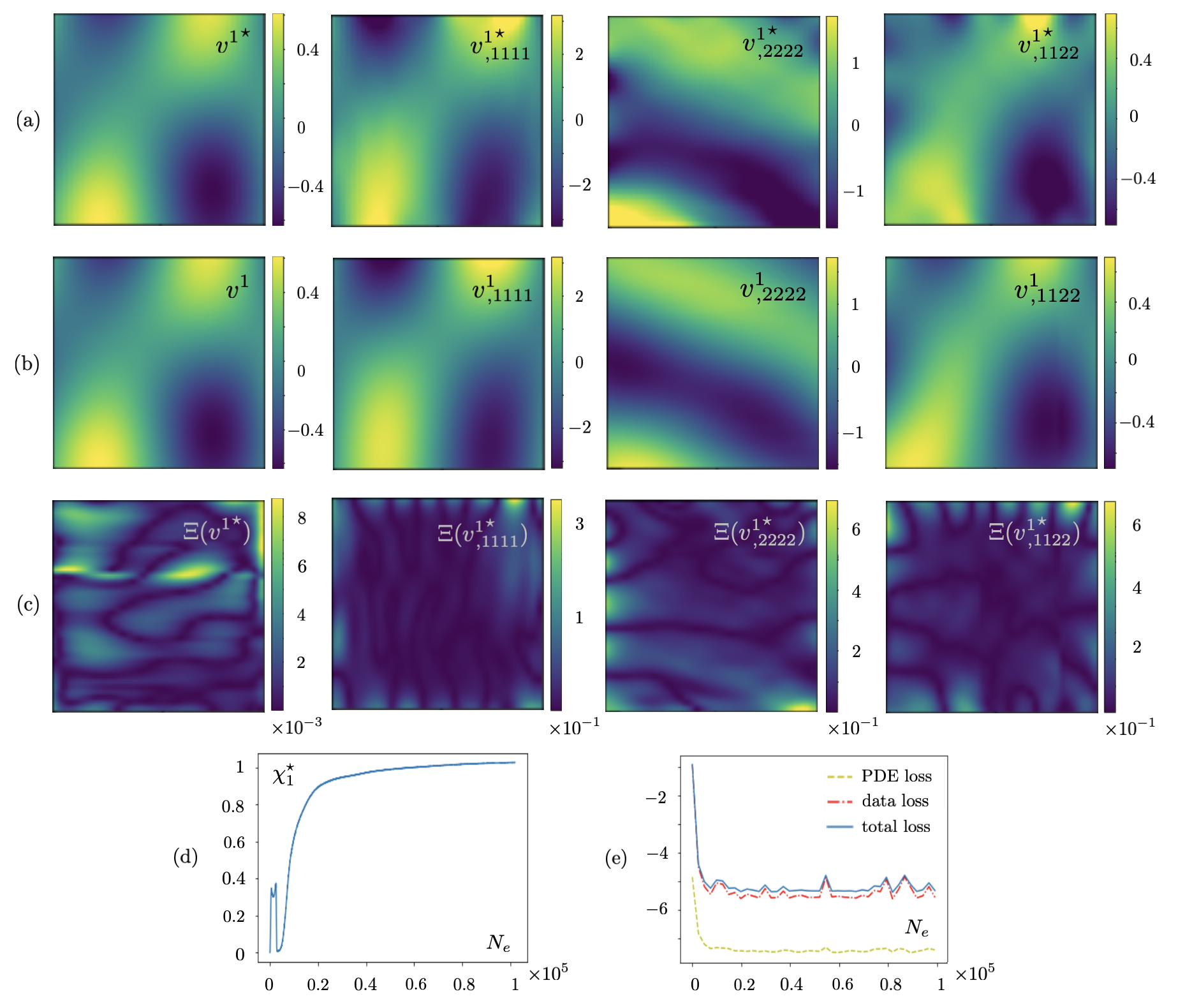} \vspace*{-3mm} 
\caption{PINN \emph{vs}.~experimental maps of particle velocity and its derivatives in $\Uppi_1^{^\text{exp}}\nxs$:~(a) PINN estimates for $\lbrace {v^1}^\star\!, {v^1}^\star_{\!\!\!\hspace{-0.1mm} ,1111}, {v^1}^\star_{\!\!\!\hspace{-0.1mm} ,2222}, {v^1}^\star_{\!\!\!\hspace{-0.1mm} ,1122} \rbrace$ wherein the derivatives are obtained by automatic differentiation, (b) normalized LDV-captured particle velocity field $v^1$ and its corresponding spectral derivatives,~(c) normal misfit~\ref{EM} between (a) and (b),~(d)~PINN-reconstructed PDE parameter $\chi^\star_1$ \emph{vs}.~the number of epochs $N_e$, and~(e) total loss $\mathcal{L}_{\varpi}$ and its components (the PDE residue and data misfit) \emph{vs.} $N_e$ in log scale. }
\label{PINNV}
\vspace*{-3mm}
\end{figure}

Fig.~\ref{PINNV} demonstrates (a) the accuracy of PINN-estimated field ${v^1}^\star$ compared to the test data $v^1$, (b) performance of automatic differentiation in capturing the fourth-order field derivatives e.g.,~${{v}^1}^\star_{\!\!,1111}$ that appear in the governing PDE~\eqref{map-exp}, and (c) the evolution of parameter $\chi_1^\star$. The comparison in (b) is with respect to the spectral derivates of test data according to Section~\ref{DIT}. It is no surprise that the automatic differentiation incurs greater errors in estimating the higher order derivatives involved in the antiplane wave motion compared to the second-order derivatives of Section~\ref{PINN-synt}.

\begin{figure}[!h]
\center\includegraphics[width=0.82\linewidth]{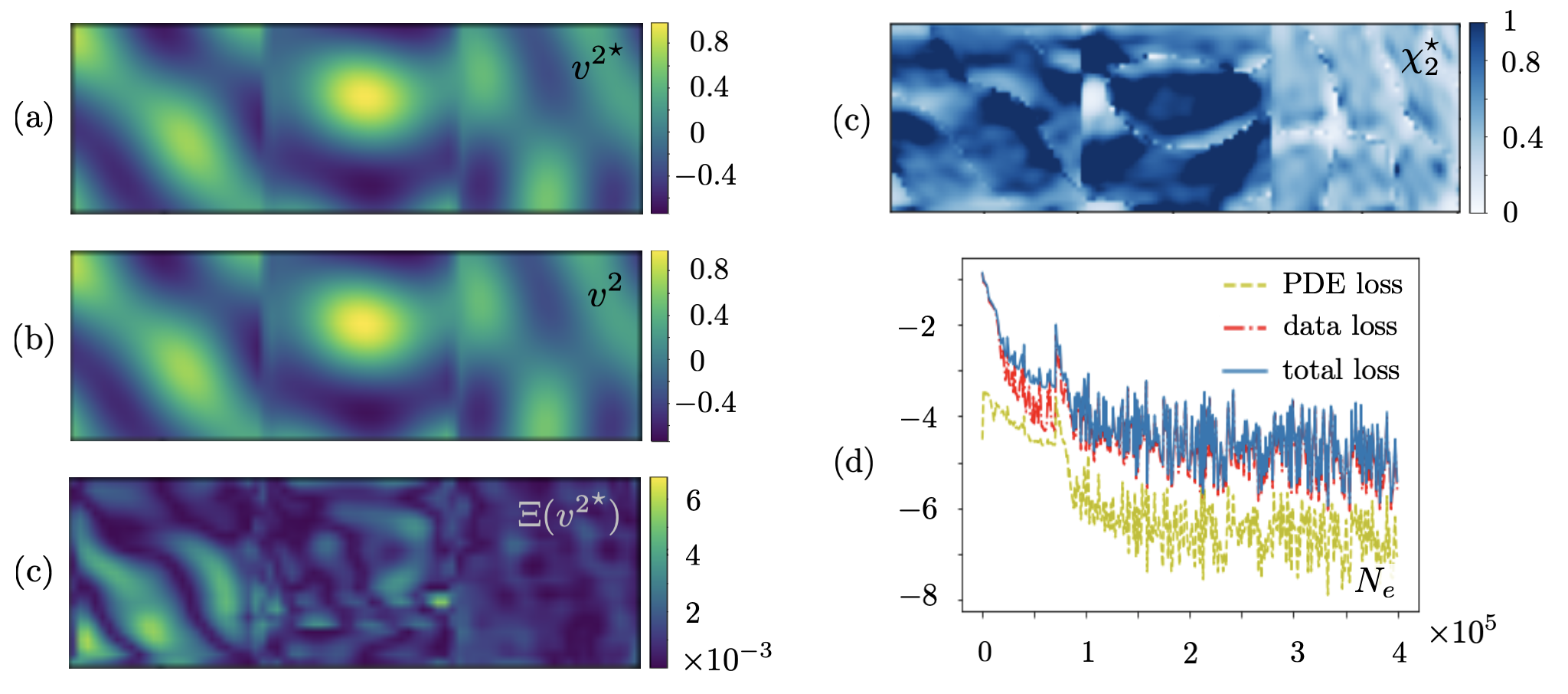} \vspace*{-3mm} 
\caption{Low-frequency PINN reconstruction in $\Uppi_2^{^\text{exp}}\!$ using test data from a single source at $f_2 = 0.17$:~(a)~PINN-predicted distribution of particle velocity ${v^2}^\star$,~(b)~normalized LDV-captured particle velocity $v^2$,~(c) normal misfit between (a) and (b),~(d)~PINN-predicted distribution of the PDE parameter $\chi^\star_2$, and~(e)~total loss $\mathcal{L}_{\varpi}$ and its components (the PDE residue and data misfit) \emph{vs.} the number of epochs $N_e$ in log scale.}
\label{PINNE}
\vspace*{-2mm}
\end{figure}

In addition, the PINN ${v^2}^\star = {v^2}^\star (\bxi,f,\bx \exs | \exs \gamma, \chi_2^\star)$ with seven hidden layers of respectively $40$, $40$, $120$, $120$, $80$, $40$, and $40$ neurons is created to reconstruct (i) particle velocity field ${v^2}$ in the layered specimen $\Uppi_2^{^\text{exp}}\nxs$, and (ii) distribution of the PDE parameter $\chi_2$ in the sampling area. The latter is defined as an unknown parameter of the network with dimension $40 \!\times\! 120$, and the scaler weight $\gamma$  is set to $\frac{1}{h (2\pi f_2)^2} = 5.84$ for the low-frequency reconstruction. In this setting, the input/output dimension for ${v^2}^\star$ reads $4800 \!\times\! 1 \!\times\! 1$. Fig.~\ref{PINNE} provides a comparative analysis between the experimental and PINN-predicted maps of velocity and PDE parameter. The associated statistics are provided in Table~\ref{Num7}. It is evident from the waveform in Fig.~\ref{PINNE}~(a) that the most pronounced errors in Fig.~\ref{PINNE}~(d) occur at the loci of vanishing particle velocity. Similar to Section~\ref{DI-synt}, this could be potentially addressed by enriching the test data. 

\section{Conclusions}

The ML-based direct inversion and physics-informed neural networks are investigated for full-field ultrasonic characterization of layered components. Direct inversion makes use of signal processing tools to directly compute the field derivatives from dense datasets furnished by laser-based ultrasonic experiments. This allows for a simplified and controlled learning process that specifically recovers the sought-for physical fields through minimizing a single-objective loss function. PINNs are by design more versatile and particularly advantageous with limited test data where waveform completion is desired (or required) for mechanical characterization. PINNs multi-objective learning from ultrasonic data may be more complex but can be accomplished via carefully gauged loss functions. 

In direct inversion, Tikhonov regularization is critical for stable reconstruction of distributed PDE parameters from limited or multi-fidelity experimental data. In this vein, deep learning offers a unique opportunity to simultaneously recover the regularization parameter as an auxiliary field which proved to be particularly insightful in inversion of experimental data.       

In training PINNs, two strategies were remarkably helpful:~(1)~identifying the reference length scale by the dominant wavelength in an effort to control the norm of spatial derivatives -- which turned out to be crucial in the case of flexural waves in thin plates with the higher order PDE, and~(2)~estimating the Lagrange multiplier by taking advantage of the inertia term in the governing PDEs.     
 
Laboratory implementations at multiple frequencies exposed that verification and validation are indispensable for predictive data-driven modeling. More specifically, both direct inversion and PINNs recover the unknown ``physical" quantities that best fit the data to specific equations (with often unspecified range of validity). This may lead to mathematically decent but physically incompatible reconstructions especially when the perceived physical laws are near their limits such that the discrepancy in capturing the actual physics is significant. In which case, the inversion algorithms try to compensate for this discrepancy by adjusting the PDE parameters which leads to non-physical reconstructions. Thus, it is paramount to conduct complementary experiments to (a) establish the applicability of prescribed PDEs, and (b) validate the predictive capabilities of the reconstructed models.     

\renewcommand{\arraystretch}{1.1}
\vspace*{-3mm}
 \begin{table}[!h]
\vspace*{-1mm}
\begin{center}
\caption{\small Mean and standard deviation of the PINN-reconstructed distributions in Fig.~\ref{PINNE} from a single-source, low-frequency test data.} \vspace*{2mm}
\label{Num7}
\begin{tabular}{|c|c|c|c|} \hline
\!$\upbeta$\! & $2_{\text{Ti}}$\! & $2_{\text{St}}$\! & $2_{\text{Br}}$\!  \\ \hline\hline  
$\chi_\upbeta$    & $0.91$  & $1$ & $0.51$  \\   
 \hline
$\exs\langle \xxs \chi_\upbeta^\star \rangle_{\Uppi_\upbeta^{^\text{exp}}\!\nxs}$     & $0.790$  & $0.890$ & $0.414$ \\ 
 \hline
$\upsigma(\chi_\upbeta^\star|_{\Uppi_\upbeta^{^\text{exp}}\!\nxs})$   & \!\!$0.214$\!\! & \!\!$0.356$\!\! & \!\!$0.134$\!\!  \\
 \hline
\end{tabular}
\end{center}
\vspace*{-3.5mm}
\end{table}

\section*{Authors' contributions} 

{\bf{Y.X.}}~investigation, methodology, data curation, software, visualization, writing -- original draft;
{\bf{F.P.}}~conceptualization, methodology, funding acquisition, supervision, writing -- original draft;
{\bf{J.S.}}~experimental data curation;
{\bf{C.W.}}~experimental data curation.


\section*{Acknowledgments} 

This study was funded by the National Science Foundation (Grant No.~1944812) and the University of Colorado Boulder through FP's startup. This work utilized resources from the University of Colorado Boulder Research Computing Group, which is supported by the National Science Foundation (awards ACI-1532235 and ACI-1532236), the University of Colorado Boulder, and Colorado State University. Special thanks are due to Kevish Napal for facilitating the use of FreeFem++ code developed as part of~\cite{pour2022(1)} for elastodynamic simulations.

\bibliography{DLels}

\end{document}